\documentclass{article}

\usepackage{amssymb,amsmath,amsfonts,amsthm}
\usepackage[letterpaper, margin=1in]{geometry}
\usepackage{graphicx}
\usepackage{url}
\usepackage{subfigure}
\usepackage{float}
\usepackage[latin1]{inputenc}
\usepackage[usenames]{color}

\usepackage{hyperref}

\newcommand{\N}{{\mathbb N}}
\newcommand{\R}{{\mathbb R}}
\newcommand{\Z}{{\mathbb Z}}
\newcommand{\zeroone}{\{ 0, 1\}}

\newcommand{\onevec}{\mathbf{1}}
\newcommand{\zerovec}{\mathbf{0}}
\newcommand{\Q}[1]{Q(#1)}
\newcommand{\Qi}[1]{Q^{*}(#1)}
\newcommand{\Qb}[1]{Q_{\beta}(#1)}

\newcommand{\PPi}[1]{P^{*}(#1)}
\newcommand{\C}[2]{C_{#1}^{#2}}
\newcommand{\ceil}[1]{\left\lceil #1 \right\rceil}
\newcommand{\floor}[1]{\left\lfloor #1 \right\rfloor}
\newcommand{\card}[1]{\left| #1 \right|}
\newcommand{\conv}{\mathrm{conv}}
\newcommand{\tabulatedset}[1]{\left\{ #1 \right\}}
\newcommand{\setof}[2]{\left\{ #1 \, : \, #2 \right\}}
\newcommand{\binomio}[2]{\left( \begin{array}{c} #1 \\ #2 \end{array} \right)}

\newcommand{\RRb}[1]{R_{\beta}(#1)}

\newcommand{\fb}{\pi_{\scriptscriptstyle \! +}}
\newcommand{\fbb}{\pi_{\scriptscriptstyle \! -}}

\newcommand{\db}{d_{\beta}}

\newcommand{\cb}{c^{+}}
\newcommand{\cbb}{c^{-}}
\newcommand{\yy}{\hat{y}}

\newcommand{\Aux}[1]{D(#1)}
\newcommand{\ab}{\bar{a}}

\newcommand{\AI}{\tilde{A}}
\newcommand{\tGam}{t(\Gamma)}

\newcommand{\tmGam}{t^{-}(\Gamma)}

\newcommand{\pGam}{p(\Gamma)}

\newcommand{\pGamp}{p(\Gamma')}

\newcommand{\Acol}[1]{{\AI}_{\bullet j}}
\newcommand{\oGam}{\circ(\Gamma)}

\newcommand{\oGamd}{\circ(\Gamma_2)}
\newcommand{\xGam}{\otimes(\Gamma)}

\newcommand{\xGami}{\otimes(\Gamma_i)}
\newcommand{\xGamu}{\otimes(\Gamma_1)}
\newcommand{\xGamd}{\otimes(\Gamma_2)}
\newcommand{\bGam}{\bullet(\Gamma)}
\newcommand{\bGami}{\bullet(\Gamma_i)}

\newcommand{\bGamd}{\bullet(\Gamma_2)}

\newcommand{\BGra}{\tilde{D}(A, \Gamma)}
\newcommand{\BVer}{\tilde{V}}
\newcommand{\BEdg}{\tilde{E}}

\def\myfnsymbol#1{\ensuremath{\ifcase#1\or \star\or \dagger\or \ddagger\or
   \mathsection\or \mathparagraph\or \|\or **\or \dagger\dagger
   \or \ddagger\ddagger \else\@ctrerr\fi}}

\newenvironment{myproof}{\begin{proof}}{\end{proof}}

\newtheorem{theorem}{Theorem}[section]
\newtheorem{corollary}[theorem]{Corollary}
\newtheorem{lemma}[theorem]{Lemma}
\newtheorem{definition}[theorem]{Definition}
\newtheorem{remark}[theorem]{Remark}
\newtheorem{assumption}[theorem]{Assumption}

\begin{document}

%

\begin{center}
         \bfseries\Large On dominating set polyhedra of circular interval graphs\footnote{Partially supported by PIP-CONICET 277, PID-UNR 416, PICT-ANPCyT 0586, and MathAmSud 15MATH06 PACK-COVER.}
 \end{center}

\smallskip

\begin{center}
\large Silvia Bianchi$^{2}$, Graciela Nasini$^{1,2}$, Paola Tolomei$^{1,2}$, and Luis Miguel Torres$^{3}$\\´
\normalsize \{\texttt{sbianchi,nasini,ptolomei}\}\texttt{@fceia.unr.edu.ar, luis.torres@epn.edu.ec} \vspace{0.5cm}\\
\small{$^{1}$CONICET - Argentina}\\
\small{$^{2}$FCEIA, Universidad Nacional de Rosario, Rosario, Argentina}\\
\small{$^{3}$Centro de Modelizaci{\'o}n Matem{\'a}tica - ModeMat, Escuela Polit{\'e}cnica Nacional, Quito, Ecuador}\\ 
\end{center}

\medskip

\begin{abstract}

Clique-node and closed neighborhood matrices of circular interval graphs are circular matrices. The
stable set polytope and the dominating set polytope on these graphs are therefore closely related to the
set packing polytope and the set covering polyhedron on circular matrices. Eisenbrand 
\emph{et al.}  \cite{EisenbrandEtAl08} take advantage of this relationship to propose a complete linear description of the stable set polytope on circular interval graphs. In this paper we follow similar ideas to obtain a
complete description of the dominating set polytope on the same class of graphs. As in the
packing case, our results are established for a larger class of covering polyhedra of the
form $\Qi{A,b}:= \conv\setof{x \in \Z^n_{+}}{Ax \geq b}$, with $A$ a circular matrix and
$b$ an integer vector.
These results also provide linear descriptions of polyhedra associated with several variants
of the dominating set problem on circular interval graphs.


\end{abstract}


\section{Introduction}

The well-known concept of domination in graphs was introduced by Berge \cite{Berge62}, modeling many facility location problems in Operations Research. Given a graph $G=(V,E)$, $N[v]$ denotes the closed neighborhood of the node $v \in V$. A set $D \subseteq V$ is called a \emph{dominating set of $G$}  if $D \cap N[v] \neq \emptyset$ holds for every $v \in V$. Given a vector $w \in \R^V$ of node weights, the \emph{Minimum-Weighted Dominating Set Problem} (MWDSP for short) consists in finding a dominating set $D$ of $G$ that minimizes $\sum_{v\in D} w_v$.
MWDSP arises in many applications, involving the strategic placement of resources on the nodes of a network. As example, consider a computer network in which one wishes to choose a smallest set of computers that are able to transmit messages to all the remaining computers \cite{KratschEtAl93}. 
Many other interesting examples include sets of representatives, school bus routing, $(r,d)$-configurations, placement of radio stations, social network theory, kernels of games, etc. \cite{HaynesEtAl98}.

The MWDSP is NP-hard for general graphs and has been extensively investigated from an algorithmic point of view (see, e.g., \cite{Bertossi84,Chang98,CorneilStewart90,Farber84}). In particular, efficient algorithms for the problem on interval and arc circular graphs are proposed in \cite{Chang98b}.

However, only a few results about the MWDSP have been established from a polyhedral point of view.
The \emph{dominating set polytope} associated with a graph $G$ is defined
as the convex hull of all incidence vectors of dominating sets in $G$.
In \cite{BouchakourEtAl08} 
the authors provide a complete description of the dominating set polytope of cycles. As a generalization, 
a description of the dominating set polytope associated with web graphs of the form
$W_{s(2k+1)+t}^{k}$, with $2\leq s \leq 3$, $0\leq t \leq s-1$, and $k \in \N$, is
presented in \cite{BianchiEtAl14b}.

Actually, the MWDSP can be regarded as a particular case of the \emph{Minimum-Weighted Set Covering Problem} (MWSCP).  Given a $\{0,1\}$-matrix $A$ of order $m\times n$, a \emph{cover} of $A$ is a vector $x\in \{0,1\}^n$ such that $Ax\geq \onevec$, where $\onevec \in \Z^m$ is the vector having all entries equal to one.
The MWSCP consists in finding a cover of minimum weight with respect to a given a weight vector $w \in \R^n$. This problem can be formulated as the integer linear program 
$$
\min\setof{w^T x}{Ax \geq \onevec, x \in \Z^n_{+}}.
$$
The set $\Qi{A}:= \conv\setof{x \in \Z^n_{+}}{Ax \geq \onevec }$ is termed as the \emph{set covering polyhedron} associated with $A$. 

The \emph{closed neighborhood matrix} of a graph $G$ is the square matrix $N[G]$ 
whose rows are the incidence vectors of the sets $N[v]$, for all $v \in V$. Observe that $x$ is the incidence vector of a dominating set of $G$ if and only if $x$ is a cover of $N[G]$. Therefore, solving the MWSCP on $N[G]$ is equivalent to solving the MWDSP on $G$. Moreover, the structure of the dominating set polytope of $G$ can be studied by considering the
set covering polyhedron associated with $N[G]$.

The closed neighborhood matrix of a web graph is a circulant matrix. More generally, the closed neighborhood matrix of a circular interval graph is a circular matrix (both terms are explained in more detail in the next section).
In this paper we are interested in studying the dominating set polytopes associated with circular interval graphs.

Another classic set optimization problem is the \emph{Set Packing Problem}: given a $\{0,1\}$-matrix $A$ of order $m\times n$, a \emph{packing} of $A$ is a vector $x \in \{0,1\}^n$ such that $Ax \leq \onevec$.
For a weight vector $w\in \mathbb R^n$, the \emph{Maximum-Weighted Set Packing Problem} (MWSPP) can
be stated as the integer linear program 
$$
\max\setof{w^Tx}{x \in \{0,1\}^n, Ax \leq \onevec}.
$$
The polytope $\PPi{A}:= \conv\setof{x \in \Z_+^n}{ Ax \leq \onevec}$ is the \emph{set packing polytope} associated with $A$.

Set packing polyhedra have been extensively studied because of their relationship with the stable set polytope. 
Indeed, given a graph $G$, a matrix $A$ can be defined whose rows
are incidence vectors of the maximal cliques in $G$. Conversely, given an arbitrary $\{0, 1\}$-matrix $A$, 
the conflict graph $G$ of $A$ is defined as a graph having one node for each column of $A$ and
two nodes joined by an edge whenever the respective columns have scalar product
distinct from zero. In both cases, stable sets in $G$ correspond to packings of $A$.

In \cite{EisenbrandEtAl08,Stauffer05} the authors present a complete linear description of the stable set polytope of circular interval graphs, 
which is equivalent to obtaining a complete linear description for the set packing polytope related to circular matrices.
The authors show that if $A$ is a circular matrix then $P^*(A)$ is completely described by three classes of inequalities: non-negativity constraints, clique inequalities, and 
\emph{clique family inequalities} introduced in \cite{Oriolo03}. 
Moreover, facet inducing clique family inequalities are associated with subwebs of the circular interval graph \cite{Stauffer05}.

Actually, their results are stated for a more general packing polyhedron $\PPi{A,b}$, defined as the convex hull of non-negative integer solutions of the system $A x \leq b$, with $b\in  \Z_{+}^m$ and $A$ a circular matrix.
In the covering case, a similar polyhedron $\Qi{A,b}$ can be defined as  the convex hull of the integer points in $\Q{A,b}=\{x\in \R^n_+ : A x \geq b\}$, with $b\in \Z_{+}^m$. When $A$ is the closed neighborhood matrix of a graph, the extreme points of $Q^*(A,b)$ correspond to some variants of dominating sets in graphs. In particular, if $b= k \mathbf{1} $, they correspond to $\{k\}$- dominating functions 
\cite{Bange} and, in the general case, they are related to $L$-dominating functions 
\cite{Lee}. Considering the symmetry in the definition of $\PPi{A,b}$ and $\Qi{A,b}$, it is natural to ask if the ideas proposed in \cite{EisenbrandEtAl08,Stauffer05} can be applied in the covering context.

In this paper we present a complete linear description of $\Qi{A,b}$ for any circular matrix $A$ and any vector $b\in \Z_n^+$. 
This yields a complete description of the polyhedron associated to $L$-dominating functions of circular interval graphs. The linear inequalities have a particular structure when $b= k \onevec$, which includes the case of $\{k\}$-dominating functions. Finally, if $k=1$, facet defining inequalities of $\Qi{A}$ provide a characterization of facets of the dominating set polytope on circular interval graphs. These inequalities are related to circulant minors of $A$.

In the light of previous results obtained by Chudnovsky and Seymour \cite{ChudnovskySeymour05},
the linear description presented in \cite{EisenbrandEtAl08,Stauffer05}
actually provided
the final piece for establishing a complete linear description of the stable set polytope for the 
much broader class of quasi-line graphs.
The fact that the dominating set problem is known to be NP-hard already for the particular subclass of line-graphs \cite{Yanakakis}, discourages seeking for an analogous result regarding domination
on quasi-line graphs. Nonetheless, we present here some positive results for a prominent subclass of them.

Some results presented in this paper appeared without proofs in \cite{Torres15}.

\section{Preliminaries}
\label{sec:preliminares}

A \emph{circular-arc graph} is the intersection graph of a set of arcs on the circle, i.e., $G=(V,E)$ is a circular-arc graph if each node $v \in V$ can be associated with an arc $C(v)$ on the circle in such a way that $uv \in E$ if and only $C(u)$ intersects $C(v)$. If additionally the family $\setof{C(v)}{v \in V}$ can be defined in such a way that no arc properly contains another,  then $G$ is a \emph{proper circular-arc graph}. Proper circular-arc graphs are also termed as \emph{circular interval graphs} in \cite{ChudnovskySeymour08} and defined in a different, but equivalent manner: 
take a finite set $V$ of points on a circle $C$ and a collection $\mathcal{I}$ of intervals from $C$.
Then, $V$ is the node set of $G$ and $u, v \in V$ are adjacent if and only if there is at least one interval in $\mathcal{I}$ containing both $u$ and $v$. Circular interval graphs are an important subclass of quasi-line graphs.
\emph{Web graphs} $W_n^k$ are regular circular interval graphs having node degree equal to $2k$.

For $n \in \N$, $[n]$ will denote the additive group defined on the set
$\tabulatedset{1, \ldots, n}$, with integer addition modulo $n$. 
Given $a,b\in[n]$, let $t$ be the minimum non-negative integer such that $a+t=b \mod n$. Then, $[a,b]_n$ denotes the \emph{circular interval} defined by the set $\{a+s: 0\leq s \leq t\}$. Similarly, $(a,b]_n$, $[a,b)_n$, and $(a,b)_n$ correspond to $[a,b]_n\setminus \{a\}$, $[a,b]_n\setminus \{b\}$, and $[a,b]_n\setminus \{a,b\}$, respectively. 

Unless otherwise stated, throughout this paper $A$ denotes a $\{0,1\}$-matrix of order $m\times n$.
Moreover, we consider the columns (resp. rows) of $A$ to be indexed by
$[n]$ (resp.~by $[m]$) and denote its entries by $a_{ij}$ with $i\in [m]$ and $j\in[n]$. 
Two matrices $A$ and $A'$ are
\emph{isomorphic}, written as $A\approx A'$, if $A'$ can be
obtained from $A$ by a permutation of rows and columns. 

In the context of this paper, a matrix $A$ is called \emph{circular} if, for every row $i \in [m]$,
there are two integer numbers $\ell_i, k_i\in [n]$ with $2\leq k_i\leq n-1$
such that the $i$-th row of $A$ is the incidence vector of the set $[\ell_i,\ell_i+k_i)_n$. 
The following is an example of a $3\times 7$-circular matrix with $\ell_1=1, \ell_2= 2, \ell_3= 5, k_1= 3, k_2=5, $ and $k_3=5$:

$$A=\left(
\begin{array}{cccccccc}
1&1&1&0&0&0&0\\
0&1&1&1&1&1&0\\
1&1&0&0&1&1&1\\
\end{array}
\right)
$$


A row $i$ of $A$ is said to \emph{dominate} a row $\ell \neq i$ of $A$ if $a_{ij} \geq a_{\ell j}$, for all
$j \in [n]$. Moreover, a row is dominating if it dominates some other row.
A square circular matrix of order $n$ without dominating rows is called
a \emph{circulant matrix}. Observe that in this case 
$k_i=k$ must hold for every row $i \in [n]$  and we can assume w.l.o.g. $\ell_i=i$ for all $i\in [n]$. 
Such a matrix will be denoted by $\C{n}{k}$. 

Given $N\subset [n]$, the \emph{minor of} $A$ \emph{obtained by
contraction of} $N$, denoted by $A/N$, is the submatrix of $A$ that
results after removing all columns with indices in $N$ and
all dominating rows.  
In this work, anytime we refer to a minor of a matrix, we mean a minor obtained by
contraction. A minor of a matrix $A$ is called a \emph{circulant minor} if it is isomorphic to 
a circulant matrix.  

Circulant minors have an interesting combinatorial characterization in terms of circuits in a particular digraph \cite{Cornuejols94}.
In fact, given a circulant matrix $\C{n}{k}$, a directed auxiliary graph $G(\C{n}{k})$ is defined by considering $n$ nodes and arcs of the form $(i,i+k)$ and $(i,i+k+1)$ for every $i\in[n]$. The authors prove that if $N\subset [n]$ induces a simple circuit in $G(\C{n}{k})$, then the matrix $\C{n}{k}/N$ is 
a circulant minor of $\C{n}{k}$. In a subsequent work, Aguilera \cite{Aguilera09} shows that $C_n^k/N$ is isomorphic to a circulant minor of $\C{n}{k}$ if and only if $N$ induces $d\geq 1$ disjoint simple circuits in  $G(\C{n}{k})$, each one having the same number of arcs of length $k$ and $k+1$.

For a matrix $A$, the \emph{fractional set covering polyhedron} is given by
$\Q{A}:= \{x \in \R^n  :  Ax \geq \onevec, \, x \geq 0\}$.
The term \emph{boolean inequality} denotes each one of the inequalities 
defining $\Q{A}$.
The \emph{covering number} $\tau(A)$ of $A$ is the minimum cardinality of a cover of $A$. When $A$ is the closed neighborhood of a graph $G$, $\tau(A)$ coincides with the \emph{domination number} $\gamma(G)$ of $G$. 
The inequality $\sum_{j=1}^n x_j\geq \tau(A)$ is called the \emph{rank constraint}, and it is always valid for $\Qi{A}$.
When $A=\C{n}{k}$ it is known that $\tau(\C{n}{k})=\left\lceil \frac{n}{k}\right\rceil$ and the rank constraint is a facet of $\Qi{\C{n}{k}}$ if and only if $n$ is not a multiple of $k$ \cite{Sassano}. In a more general sense, given a matrix $A$ and a vector $b\in \Z^m$, we define $\tau_b(A):=\min \{\onevec^T x : x\in \Qi{A,b}\}$ and the rank constraint $\sum_{j=1}^n x_j\geq \tau_b(A)$, which is always valid for $\Qi{A,b}$. When $A$ is the closed neighborhood matrix of a graph $G$ and $b=k\onevec$, $\tau_b(A)=\gamma_{\{k\}}(G)$ is the $\{k\}$-\emph{domination number} of $G$. For general $b\in \Z_n^+$, $\tau_b(A)$ is the $L$-\emph{domination number} of $G$, for the corresponding list $L$ associated to $b$.


The class of \emph{row family inequalities} (rfi) was proposed in \cite{ArgiroffoBianchi10} 
as a counterpart to clique family inequalities in the set packing case \cite{Oriolo03}.
We describe them at next, slightly modified to fit in our
current notation.

Let 
$F\subset [m]$ be a set of row indices of $A$, $s:=\card{F}\geq 2$, 
$p \in [s-1]$ such that $s$ is not a multiple of $p$, and $r:= s-p \left\lfloor \frac{s}{p}\right\rfloor$. 
Define the sets 
$$
I(F, p) = \setof{j \in [n]}{\sum_{i \in F} a_{ij} \leq p}, \quad
O(F, p) = \setof{j \in [n]}{\sum_{i \in F} a_{ij} = p+1}.
$$

 Then, the \emph{row family inequality} (rfi) \emph{induced by} $(F, p)$ is 
\begin{equation} \label{eq:rfi}
(r+1) \!\!\! \!\!\sum_{j \in O(F, p)} \!\!\! x_j  \, + \, r \!\!\!\!\! \sum_{j \in I(F, p)} \!\!\! x_j  \geq r \ceil{\frac{s}{p}}.
\end{equation}

Row family inequalities generalize several previously known classes of valid inequalities
for $\Qi{A}$. However, in contrast to clique family inequalities, not all of them are 
valid for $\Qi{A}$. In \cite{ArgiroffoBianchi10} it is proved that inequality \eqref{eq:rfi} is valid for $\Qi{A}$ if the following condition holds for every cover
$B$ of $A$:
\begin{equation} \label{eq:rfi-condition}
p \card{B \cap I(F, p)} + (p+1) \card{B \cap O(F, p)} \geq s.
\end{equation}

In particular, if $p^*:= \max_{j\in [n]} \sum_{i\in F} a_{ij}-1$, the row family inequality induced by $(F, p^*)$ is always valid for $Q^*(A)$. Throughout this article, we are going to refer to this inequality simply as the
\emph{row family inequality induced by $F$}.

%
%



In the particular case when $A=\C{n}{k}$, facet defining inequalities of $\Qi{C^k_n}$ related to circulant minors were studied in \cite{ArgiroffoBianchi09,BianchiEtAl14a,BianchiEtAl14b,TolomeiTorres15,Torres15}. Given $N \subset [n]$ such that 
$\C{n}{k}/ N \approx \C{n'}{k'}$,
let $W:=\{j\in N : j-(k+1) \in N\}$.
Then, the  inequality 
\begin{equation}
\label{eq:minor-eq}
2 \sum_{j\in W} x_j + \sum_{j\notin W} x_j  \geq \ceil{\frac{n'}{k'}}
\end{equation}
is valid for $\Qi{\C{n}{k}}$, and facet defining if $n'-k' \left\lfloor \frac{n'}{k'}\right\rfloor=1$. 
This inequality is termed as the \emph{minor inequality} induced by $N$ \cite{ArgiroffoBianchi09,BianchiEtAl14a}. 

For a general circular matrix $A$, if $A/ N \approx \C{n'}{k'}$ and $F \subset [m]$ is the set of
rows of $A/ N$, then the rfi induced by $F$ will be termed as \emph{minor related row family inequality}.
These inequalities were introduced in \cite{Torres15} for the specific case when $A$ is a circulant
matrix. In this setting, the inequalities can be seen as a generalization of the minor inequalities \eqref{eq:minor-eq}, as they have the form:
\begin{equation}
\label{eq:rfi-minor}
(r+1) \sum_{j\in W} x_j + r \sum_{j\notin W} x_j \geq r \left\lceil \frac{n'}{k'}\right\rceil.
\end{equation}
with $r= n'-k' \left\lfloor \frac{n'}{k'}\right\rfloor$. 


In this paper we follow many 
of the ideas proposed in \cite{EisenbrandEtAl08,Stauffer05} for describing the stable set polytope of circular interval graphs. Actually, the construction detailed below was originally presented by
Bartholdi, Orlin and Ratliff \cite{BartholdiEtAl80} in the context of an algorithm to solve the cyclic
staffing problem, which is equivalent to the task of minimizing a linear function over $\Qi{A,b}$.

 We associate with a circular matrix the digraph defined as follows: 
\begin{definition}\label{D(A)}
Given a 
circular  matrix $A$, let $\Aux{A}:=(V, E)$ where 
$V:=[n]$ and $E$ 
is the union of the following four sets:
\begin{align*}
E_1^+&:= 
\setof{a_i:= (\ell_i -1, \ell_i+k_i-1)}{ i\in [m]}, \\
E_2^{+} &:= 
\setof{a_{m+j}:= (j -1, j)}{ j \in [n]}, \\
E_1^{-} &:= 
\setof{\ab_i:= (\ell_i+k_i-1, \ell_i -1)}{ i \in [m]}, \\
E_2^{-} &:= 
\setof{\ab_{m+j}:= (j, j-1)}{ j \in [n]}.
\end{align*}
The arcs in
$E_1^{+} \cup E_2^{+}$ 
are called \emph{forward arcs}, while the arcs in  $E_1^{-} \cup E_2^{-}$
are \emph{reverse arcs}. Moreover, arcs in $E_1^{+} \cup E_1^{-}$ are termed as \emph{row arcs}, and
arcs in $E_2^{+} \cup E_2^{-}$ are \emph{short arcs}.

For any path $P$ in $D(A)$, $E(P)$ denotes the set of arcs from $P$, whereas $E^+_1(P)$, $E^+_2(P)$, $E^-_1(P)$, and $E^-_2(P)$ denote the sets $E^+_1\cap E(P)$, $E^+_2\cap E(P)$, $E^-_1\cap E(P)$, an $E^-_2\cap E(P)$, respectively.

The (oriented) \emph{length} $l(a_i)$ (resp.~$l(\ab_i)$)
of an arc $a_i \in E_1^{+}$ (resp.~$\ab_i \in E_1^{-}$) is equal to $k_i$ (resp.~to $-k_i$). Arcs in $E_2^{+}$
have length of 1, while the length of arcs in $E_2^{-}$ is equal to -1.

\end{definition}

Simple directed circuits in $D(A)$ play a important role in the description of the set packing polytope associated with $A$ and the more general packing polytope $P^*(A,b)$
defined in the last section \cite{EisenbrandEtAl08,Stauffer05}. In this paper 
we show that similar results hold for the corresponding covering polyhedra $\Qi{A}$ and $\Qi{A,b}$.
Throughout this article we will use the term \emph{circuit} to refer to a simple directed circuit.

Consider the invertible linear map $T: \R^n \to \R^n$ represented by a $\{0, 1, -1\}$-matrix $T$ having the elements on the diagonal all equal to 1, the elements on the first subdiagonal equal to -1 and all other elements equal to zero, i.e.,
\begin{equation}
\label{matrizT}
T := \left(
\begin{array}{cccc}
1 &&&\\
-1 & 1 &&\\
& \ddots & \ddots &  \\
& & -1 & 1 
\end{array}
\right).
\end{equation}
For a circular matrix $A$,  let $\AI:=\binomio{A}{I} \in \{0,1\}^{(m + n) \times n}$, with $I$ being the identity matrix of order $n$.
If $B:= \AI T$ and
$M$ denotes the submatrix consisting of the first $n-1$ columns of $B$, it  
is straightforward to verify that the node-arc incidence matrix $H$ of the digraph $\Aux{A}$ is given by
\begin{equation}
\label{eq:node-arc-incidence}
H:=
\left(
\arraycolsep=1.4pt\def\arraystretch{2.2}
\begin{array}{r|r}
M^T & - M^T \\[0.3\baselineskip]
\hline
-\onevec^T M^T & \onevec^T M^T
\end{array}
\right).
\end{equation}



The remaining of this article is structured as follows. In the next section, we 
review some constructions and results presented in \cite{EisenbrandEtAl08,Stauffer05} in the context of the covering case. 
In Section~\ref{sec:gamma-ineq-suff}
we introduce a class of valid inequalities for $\Qi{A,b}$ induced by circuits in $D(A)$
and show that these inequalities are sufficient for describing $\Qi{A,b}$, for any $b \in \Z_{+}^m$. 
From this result we obtain a complete description of the polyhedron of $L$-dominating
functions on circular interval graphs.
In Section~\ref{sec:homog-rhs}
we consider polyhedra of the form $\Qi{A,k \onevec}$ with $k \in \N$, which include the set covering polyhedron as a particular case. We prove that, for this class of polyhedra,
the circuits in $D(A)$ inducing facet defining inequalities have no reverse row arcs.
In Section~\ref{sec:no-reverse-arcs} we further study the structure of such circuits, obtaining as a result that the corresponding inequalities have full support and only two consecutive positive integer coefficients. 
In the particular case of the set covering problem, these inequalities are row family inequalities.
Finally, in Section~\ref{sec:set-covering-polyhedron} we prove that the relevant inequalities are related to circulant minors.  As we have observed in the introduction, the description of $\Qi{A,k \onevec}$
yields a complete description of the polyhedra associated with $\{k\}$-dominating functions on circular interval graphs, whereas a complete description for the dominating set polytope on those graphs can be obtained from the description of $\Qi{A}$.


%
%
%
%
%

\section{Following the ideas of the packing case}


As we have mentioned, the study of the covering polyhedra of circular matrices closely follows the ideas proposed in \cite{EisenbrandEtAl08,Stauffer05} for the corresponding packing polytopes. Some of these ideas can be straightforwardly translated to the covering case. In this section we review them, including the corresponding proofs for the sake of completeness.   

Given a circular matrix $A$, $b\in \Z_+^m$ and $\beta \in \N$, the \emph{slice of} $Q(A,b)$ \emph{defined by} $\beta$ is the polyhedron:
$$
\Qb{A,b}:= \Q{A,b} \cap \setof{x \in \R^n}{\onevec^T x =\beta}.
$$
Remind that $\AI=\binomio{A}{I}$, $B=\AI T$, and $M$ is the submatrix consisting of the first $n-1$ columns of $B$. Moreover let $d=\binomio{b}{\mathbf{0}}$ and let $v$ be the last column of $\AI$.

\begin{lemma}
\label{th:int-slices}
For any circular matrix $A$, $b\in \Z_+^m$ and $\beta \in \N$, the polytope $\Qb{A,b}$
is integral.
\end{lemma}

\begin{myproof}

Let $\RRb{A,b}$ be the image of $\Qb{A,b}$ under the inverse of matrix $T$ defined in \eqref{matrizT}, i.e.,
$\RRb{A,b}=\setof{T^{-1}x}{x\in \Qb{A,b}}$. Then if $y:= T^{-1}x$ we have 
\begin{align*}
\RRb{A,b} &= \setof{y \in \R^n}{By \geq d, y_n = \beta} \\
&= \setof{y \in \R^n}{M \yy + \beta v \geq d} \\
&= \setof{(\yy,\beta) \in \R^n}{M \yy \geq \db},
\end{align*}
where $\db:= d - \beta v \in \Z^{m+n}$ and $\yy \in \R^{n-1}$ denotes the vector obtained from $y$ by dropping its last coordinate $y_n$. 
Since $M^T$ is a submatrix of the node-arc incidence matrix of digraph $D(A)$, it follows that $M$ is totally unimodular. Thus, $\RRb{A,b}$ is integral. Moreover, since $T$ maps integral points onto integral points, it follows that $\Qb{A,b}$ is also integral.
\end{myproof}

\begin{corollary}
\label{th:int-slices-cor}
If $A$ is a circular matrix  and $b\in \Z_+^m$ then
$$
\Qi{A,b}= \conv \left( \bigcup\limits_{\beta \in \N} \Qb{A,b} \right).
$$
\end{corollary}

The last result states that the \emph{split-rank} of the polyhedron $\Q{A,b}$ is equal to
one. This fact can be used to address the problem of separating a point from $\Qi{A,b}$. To do that, we need the following definitions:

\begin{definition}
Let $A$ be a circular matrix, 
$b\in \Z_+^m$ and $x^{*} \in \Q{A,b}$.  Consider the \emph{slack vector} $s^{*}:= \AI x^{*} -d \geq 0$ and let $\mu:= \ceil{\onevec^T x^*} - \onevec^T x^*$. We define the cost vectors $\cb(x^*), \cbb(x^*) \in \R^{m+n}$ by
\begin{equation}
\label{eq:arc-costs}
\begin{aligned}
\cb(x^*) &:= \mu (s^{*} - (1 - \mu) v),\\ 
\cbb(x^*) &:= (1 -\mu) (s^{*} + \mu v). 
\end{aligned} 
\end{equation}
\end{definition}

From now on, $D(A,x^*)$ denotes the digraph $D(A)$ with cost $c^+(x^*)$ on its forward arcs and cost $c^-(x^*)$ on its reverse arcs. For any path $P$ in $D(A)$, $c(P,x^*)$ denotes the cost of $P$ in $D(A,x^*)$.

\begin{remark}\label{noneg}
Observe that if $\onevec^T x^* \in \Z$ then $\cb(x^*)=0$, $\cbb(x^*)= s^{*} \geq 0$ and 
then $D(A, x^*)$ cannot contain a negative cost circuit.
\end{remark}

Non-negativity of all circuits in $D(A,x^*)$ is a sufficient condition for a point $x^*\in \Q{A,b}$ to be in $\Qi{A,b}$, as shown at next.

\begin{lemma} \label{th:gammanonegativo}
Let $A$ be a circular matrix, 
$b\in \Z_+^m$ and $x^* \in \Q{A,b}$ such that every circuit in $D(A,x^*)$ has non-negative cost.
Then $x^* \in \Qi{A,b}$. 
\end{lemma}

\begin{myproof}
From Corollary~\ref{th:int-slices-cor} it is sufficient to consider the case $\onevec^T x^* \notin \Z$.
Let $\tau^*:= \floor{\onevec^T x^*}$ and $\mu:= \tau^* + 1 - \onevec^T x^*>0$.

Since $D(A,x^*)$ does not contain a negative cost circuit, there exists a vector of \emph{node potentials} $z \in \R^n$,
such that $z_i - z_{j}$ is at most the cost of arc $(i,j)$, 
for every arc $(i,j) \in E$ (see, e.g. \cite[Chapter 2]{CookEtAl98} for a proof of this well-known result).  
Considering the node-arc incidence matrix $H$ defined in \eqref{eq:node-arc-incidence}, this property  
can be written as 
$$H^T z \leq \binomio{\cb(x^*)}{\cbb(x^*)}$$ 
or equivalently, 
\begin{equation}
\label{eq:potential-z}
-c^-(x^*)\leq  M \hat{z} - z_n M \onevec \leq c^+(x^*), \\
\end{equation}
where $\hat{z}$ denotes the vector obtained from $z$ by dropping its last coordinate $z_n$. 

Now define
$$
\begin{aligned}
x^1&:= x^* - \frac{1}{\mu}Tz + \frac{z_n}{\mu} e_1 - (1 - \mu) e_n, \\
x^2&:= x^* + \frac{1}{1-\mu}Tz - \frac{z_n}{1-\mu} e_1 + \mu e_n. 
\end{aligned}
$$ 
It is straightforward to verify that $x^* = \mu x^1 + (1- \mu) x^2$. Thus, if $x^1, x^2 \in \Qi{A,b}$ then
$x^* \in \Qi{A,b}$ follows from convexity of this polyhedron. Moreover, since $\onevec^T T z = z_n$, we have
$\onevec^T x^1 = \onevec^T x^* - (1 - \mu) = \tau^*\in \N$ and $\onevec^T x^2 = \onevec^T x^* + \mu = \tau^* + 1\in \N$.
Hence, due to Corollary~\ref{th:int-slices-cor}, it suffices to show that $x^1, x^2 \in \Q{A,b}$. Indeed,
\begin{align*}
\AI x^1 &= \AI x^* - \frac{1}{\mu}(\AI Tz - z_n \AI e_1) - (1 - \mu) \AI e_n, \\
&= \AI x^* - \frac{1}{\mu}\left(M\hat{z} - z_n(\AI e_1 - v) \right) - (1 - \mu) v, \\
&= \AI x^* - \frac{1}{\mu}(M\hat{z} - z_n M \onevec) - (1 - \mu) v, \\
& \geq \AI x^* - \frac{1}{\mu}c^+(x^*) - (1 - \mu) v= \AI x^* - s^*= d
\end{align*}
where the third equality follows from the fact that 
$$M \onevec = B (\onevec - e_n) = \AI T (\onevec - e_n) = \AI T \onevec - v = \AI e_1 - v,$$
and the inequality in the fourth row is obtained from \eqref{eq:potential-z}. With a similar argument, $\AI x^2 \geq d$ follows, and the proof is completed. 
\end{myproof}

As an immediate consequence of the previous result, if $x^*\in\Q{A,b}\setminus\Qi{A,b}$ then there exists a negative cost circuit in $D(A,x^*)$. It follows that, similarly as observed in \cite{EisenbrandEtAl08} in the packing context, the membership problem for $\Qi{A,b}$ can be reduced to a minimum cost circulation problem in $D(A)$. 

In the next section, we present valid inequalities for $\Qi{A,b}$ associated with circuits in $D(A)$. We will see that, given $x^*\in \Q{A,b}$, if $c(\Gamma, x^*)<0$ holds for some circuit $\Gamma$ in $D(A,x^*)$, then there is a separating split cut for $x^*$ associated with $\Gamma$. In this way, we will prove that $\Qi{A,b}$ can be described by these \emph{circuit inequalities} together with the inequalities defining $\Q{A,b}$.

\section{A complete linear description of $\Qi{A,b}$}
\label{sec:gamma-ineq-suff}



%
Consider a circular matrix $A$ and the directed graph $D(A)$ defined in Section \ref{sec:preliminares}. Recall from Definition \ref{D(A)} that for any arc $a$ in $D(A)$, $l(a)$ denotes its (oriented) length.  

Given a closed directed (not necessarily simple) path $\Gamma=(V(\Gamma),E(\Gamma))$ in $D(A)$, its 
\emph{winding number} is the integer 
$\pGam$ such that:
$$
\pGam \, n = \sum_{a \in E(\Gamma)} \!\! l(a).
$$
For every $i\in [m]$, let $P_i^+$ (resp. $P_i^-$) be the path of short forward (resp. reverse) arcs in $D(A)$ connecting $l_i-1$ with $l_i+k_i-1$ (resp. $l_i+k_i-1$ with $l_i-1$).

We say that a forward row arc $a_i=(\ell_i -1, \ell_i+k_i-1) \in E_1^{+}$ \emph{jumps over} a node $j \in V$
if $j \in 
[l_i, l_i+k_i)_n$ and the only forward short arc jumping over $j$ is the arc $(j-1,j)$. A reverse arc is said
to jump over a node $j \in [n]$ if and only if the corresponding (antiparallel) forward arc jumps over $j$.

Let $p^{+}(\Gamma, j)$ and $p^{-}(\Gamma, j)$ the number of forward and reverse arcs of $\Gamma$ jumping
over a node $j \in [n]$, respectively. We have the following result:

\begin{lemma}
\label{th:winding-number}
Let $A$ be a circular matrix and $\Gamma$ be a closed directed path in $D(A)$.
For any $j \in [n]$, $$p^{+}(\Gamma, j) - p^{-}(\Gamma, j)=\pGam.$$
\end{lemma} 

\begin{myproof}
Let us start with the case when $\Gamma$ has only short arcs. For any $j\in [n]$, the arcs in $\Gamma$ that may leave $j$ are $(j,j+1)$ and $(j,j-1)$. Then the number of arcs in $\Gamma$ leaving $j$ is $p^{+}(\Gamma, j+1)+ p^{-}(\Gamma, j)$. Similarly, the number of arcs in $\Gamma$ entering $j$ is $p^{+}(\Gamma, j)+ p^{-}(\Gamma, j+1)$. Since $\Gamma$ is a closed path, we have:
$$p^{+}(\Gamma, j+1)+ p^{-}(\Gamma, j)= p^{+}(\Gamma, j)+ p^{-}(\Gamma, j+1)$$
or, equivalently, 
$$p^{+}(\Gamma, j)- p^{-}(\Gamma, j)= p^{+}(\Gamma, j+1) - p^{-}(\Gamma, j+1).$$

Hence, $\gamma= p^{+}(\Gamma, j)- p^{-}(\Gamma, j)$ is a fixed value for all $j\in [n]$.

For each $j\in [n]$, let $\delta^+(j)$ be the set of arcs of $\Gamma$ leaving $j$. We have:

\begin{align*}
\pGam \, n &= \sum_{a \in E(\Gamma)} \!\! l(a)= \sum_{j\in [n]} \sum_{a\in \delta^+(j)} \!\! l(a)= \sum_{j\in [n]} \!\! [p^{+}(\Gamma, j+1)-p^{-}(\Gamma, j)]\\ 
&= \sum_{j\in [n]} \!\! p^{+}(\Gamma, j+1)- \sum_{j\in [n]} \!\! p^{-}(\Gamma, j)= 
\sum_{j\in [n]} \!\! [p^{+}(\Gamma, j)-p^{-}(\Gamma, j)]= n \gamma.
\end{align*}  

Then, $\pGam = \gamma= p^{+}(\Gamma, j)-p^{-}(\Gamma, j)$ for all $j\in [n]$.

%
%
%

Now consider any closed path $\Gamma$ in $D(A)$ and let $\Gamma'$ be the path obtained from $\Gamma$ by replacing each forward row arc $a_i= (\ell_i -1, \ell_i+(k_i-1))$
(resp. reverse row arc $\ab_i= (\ell_i+(k_i-1), \ell_i -1)$)  
by the path $P_i^+$ (resp. $P_i^-$). Observe that $\pGam= \pGamp$.
Moreover, for each $j\in [n]$, as each row arc jumping over $j$ is replaced by a path containing exactly one short arc jumping over $j$, we have $p^{+}(\Gamma, j) = p^{+}(\Gamma', j)$ and $p^{-}(\Gamma, j) = p^{-}(\Gamma', j)$.
This completes the proof.
\end{myproof}

Given a closed directed path $\Gamma$ in $D(A)$, 
denote by $\fb \in \Z_+^{m+n}$ (resp. $\fbb\in \Z_+^{m+n}$) the vector  whose components are the number of times each forward (resp. reverse) arc in $D(A)$ occurs in $\Gamma$. In particular, if $\Gamma$ is a circuit then $\fb, \fbb \in \zeroone^{m+n}$. Observe that, for any $j\in [n]$, if $e_j$ is the $j$-th canonical vector, $p^{+}(\Gamma, j) = \fb^T \AI e_j$ and $p^{-}(\Gamma, j) = \fbb^T \AI e_j$, where $\AI=\binomio{A}{I}$. Hence, as a consequence of Lemma \ref{th:winding-number}, we have: 

\begin{corollary}
\label{th:parameter-properties}
Let $A$ be a circular matrix. Then, for any circuit $\Gamma$ in $D(A)$,
$$
(\fb - \fbb)^T \AI =  \pGam \onevec^T.
$$
\end{corollary}

As a consequence of Lemma \ref{th:gammanonegativo}, for any $x^*\in \Q{A,b}\setminus \Qi{A,b}$, there exists a negative cost circuit in $D(A,x^*)$. As we see at next, this circuit has positive winding number.



\begin{lemma}\label{ppositivo}
Let $A$ be a circular matrix, $b\in \Z_m^+$, $x^*\in \Q{A,b}\setminus \Qi{A,b}$, and $\Gamma$ a circuit with negative cost in $D(A,x^*)$. Then $p(\Gamma)>0$.
\end{lemma}
\begin{myproof}
Let $p=p(\Gamma)$. From definition we have that  $c(\Gamma, x^*) = \fb ^T\cb(x^*) + \fbb^T \cbb (x^*)$.  Since $\cb(x^*) + \cbb (x^*)=s^*$, we have
$$c(\Gamma, x^*) = \fbb^T s^* + (\fb^T-\fbb^T) \cb (x^*).$$
In addition 
$$(\fb^T-\fbb^T) \cb (x^*)= (\fb^T-\fbb^T)(\mu s^*-\mu (1-\mu)v).$$
Since $(\fb^T-\fbb^T)v=p$ holds from Corollary \ref{th:parameter-properties}, we have that
$$c(\Gamma, x^*) = \fbb^T s^* + \mu (\fb^T-\fbb^T)s^* -\mu(1-\mu)p$$
or, equivalently,
$$c(\Gamma, x^*) = (1-\mu)\fbb^T s^* + \mu \fb^T s^* -\mu(1-\mu)p.$$
Observe that $\fbb^T s^*\geq 0$, $\fb^T s^*\geq 0$ and $0<\mu<1$. It follows that  
$$c(\Gamma, x^*) \geq  -\mu(1-\mu)p.$$
If $p\leq 0$ then $c(\Gamma, x^*)\geq 0$, contradicting the assumption. Therefore, $p>0$.
\end{myproof}


Given a circuit $\Gamma$ in $D(A)$, $b\in \mathbb Z_m^+$, and $d=\binomio{b}{0}$, we introduce the following parameters:
\begin{align*}
t(\Gamma,b)&:= \sum_{i\in E^+_1(\Gamma)} \!\!\! b_i - \sum_{i\in E^-_1(\Gamma)}\!\!\!b_i= (\fb - \fbb)^T d,\\ 
\beta(\Gamma,b)&:= \floor{\frac{t(\Gamma,b)}{\pGam}}, \text{ and }\\
r(\Gamma,b)&:= t(\Gamma,b) - \beta(\Gamma,b) p(\Gamma).
\end{align*}

The parameters above are involved in the next definition:

\begin{definition}
\label{gammaineq}
Let $A$ be a circular matrix, $b\in \Z^m_+$ and $\Gamma$ be a circuit in $D(A)$. If 
$\beta=\beta(\Gamma,b)$ and $r=r(\Gamma,b)$, 
the $\Gamma$-inequality is defined as 
\begin{equation}
\label{eq:gamma-ineq}
\sum_{j\in [n]} [p^-(\Gamma,j)+r] \, x_j \geq  r \, (\beta +1)+ \sum_{i\in E^-_1(\Gamma)}\!\!\!b_i.
\end{equation}

We say that an inequality is a circuit inequality if it is a $\Gamma$-inequality for some circuit $\Gamma$ in $D(A)$.
\end{definition}


The next result shows that every $\Gamma$-\emph{inequality} with $\pGam > 0$ is valid for the slices $Q_{\beta}(A,b)$ and $Q_{\beta+1}(A,b)$ for $\beta=\beta(\Gamma,b)$ i.e., it is a disjunctive cut for $\Q{A,b}$ and then, it is valid for $\Qi{A,b}$.
\begin{theorem}\label{validas}
Let $A$ be a circular matrix, $b\in \mathbb Z_m^+$ and  $\Gamma$ be a circuit in $D(A)$ with $\pGam > 0$. Then, the $\Gamma$-\emph{inequality} 
is valid for $\Qi{A,b}$.
\end{theorem}

\begin{myproof}
In the following we denote $\pGam$, $t(\Gamma, b)$, $\beta(\Gamma,b)$, and $r(\Gamma,b)$ simply by $p$, $t$, $\beta$, and $r$, respectively. Remind that $t= (\fb - \fbb)^T d$ and let $t^-= \sum_{i\in E^-(\Gamma)}b_i= \fbb^T d$. Moreover, observe that $p>0$ implies $r\geq 0$ and $t<p (\beta+1)$.

Let $x^*$ be an extreme point of $\Qi{A,b}$. Applying $\fbb$ as a vector of  multipliers on the system $\AI x^* \geq d$, we obtain  
\begin{equation} \label{ineq}
\fbb^T \AI x^* \geq  \fbb^T d.
\end{equation}

Since $\beta=\left\lfloor \frac{t}{p}\right\rfloor \in \Z$, $x^*$ satisfies the disjunction
$$
\onevec^T x^* \geq \beta + 1 \qquad
\text{ or } \qquad
\onevec^T x^* \leq \beta.
$$
Assume at first that $\onevec^T x^* \geq \beta + 1$. 
Multiplying this inequality  by the non-negative factor $r$ and adding it with inequality (\ref{ineq}) yields

\begin{equation}
\label{eq:proof-gamma-ineq1}
\fbb^T \AI x^* + r \onevec^T x^* \geq \fbb^T d + (\beta + 1)r = t^- + (\beta +1)r.
\end{equation}

By using the fact that $\fbb^T \AI e_j=p^{-}(\Gamma, j)$ for every $j\in[n]$, we conclude that
$x^*$ satisfies the $\Gamma$-inequality (\ref{eq:gamma-ineq}).

Now suppose $\onevec^T x^* \leq \beta$. Multiplying this inequality by the negative factor $t - (\beta + 1) p$ 
and adding the valid inequality $\fb^T \AI x^* \geq \fb^T d$, we
obtain

\begin{equation}
\label{eq:proof-gamma-ineq}
\fb^T \AI x^* + [t - (\beta+1) p] \onevec^T x^* \geq \fb^T d + \beta [t - (\beta+1) p].
\end{equation}

By Corollary~\ref{th:parameter-properties} we have that  $\fb^T \AI x^*= \fbb^T \AI x^* +  p \onevec^T x^*$ and then the left-hand side of this inequality is:
$$\fb^T \AI x^* + [t - (\beta+1) p] \onevec^T x^* =
\fbb^T \AI x^* + r \onevec^T x^*.$$
Moreover, as $\fb^T d = t+t^-$  the right-hand side of \eqref{eq:proof-gamma-ineq} can be written as:
\begin{align*}
\fb^T d + \beta [t - (\beta+1) p]
& = t+t^-+ \beta t - \beta(\beta +1) p \\
& = t^- + (\beta + 1)r.
\end{align*}
Hence, $x^*$ does also fulfill \eqref{eq:gamma-ineq} when $\onevec^T x^* \leq \beta$.
\end{myproof}

The following lemma establishes a necessary condition for a circuit $\Gamma$ so that the $\Gamma$-inequality defines a facet of $\Qi{A,b}$.

\begin{lemma} \label{rem:non-redundant}
Let $\Gamma$ be a circuit in $D(A)$ such that the $\Gamma$-inequality induces a facet of $\Qi{A,b}$. Then,  
\begin{equation}
\label{eq:cond-gamma}
\pGam \mbox{ does not divide } 
t(\Gamma,b) \text{ and }  2 \leq \pGam \leq 
t(\Gamma,b)-1.
\end{equation}
Then, every circuit inequality defining a facet of $\Qi{A,b}$ has full support (i.e., non zero coefficients for all variables). 

\end{lemma}

\begin{myproof}
Note that if $\pGam$ divides $t(\Gamma,b)$ then $r(\Gamma, b)=0$ and the $\Gamma$-inequality \eqref{eq:gamma-ineq} reduces to  
$\fbb^T \AI x \geq \tmGam$, which is redundant since it is $\AI x \geq d$ multiplied by the vector $\fbb^T$. It follows that  $\pGam \geq 2$.

Moreover, if 
$\beta \leq 0$, the inequality $\onevec^T x \geq \beta +1$ is implied by any of the inequalities in the system $Ax \geq b$ and the $\Gamma$-inequality is valid for $\Q{A,b}$. Thus, $\beta \geq 1$ and since $p(\Gamma)\geq 2$ we obtain 
$\pGam \leq t(\Gamma,b)-1$.
\end{myproof}

Observe that if $x^*\in \Q{A,b}\setminus \Qi{A,b}$ and $\Gamma$ is a circuit with negative cost in $D(A,x^*)$, then from Lemma \ref{ppositivo} and Theorem \ref{validas} the $\Gamma$-inequality is valid for $\Qi{A,b}$.
We see at next that $x^*$ violates this inequality.

\begin{lemma}
\label{5prima}
Let $A$ be a circular matrix, $b\in \Z_m^+$, $x^*\in \Q{A,b}\setminus \Qi{A,b}$ and $\Gamma$ be a circuit with negative cost in $D(A,x^*)$. Then, the $\Gamma$-inequality is 
violated by $x^*$.  
\end{lemma}

\begin{myproof}
In the following we denote $\pGam$, $t(\Gamma, b)$, $\beta(\Gamma,b)$, and $r(\Gamma,b)$ simply by $p$, $t$, $\beta$, and $r$, respectively. Remind that $t= (\fb - \fbb)^T d$ and let $t^-= \sum_{i\in E^-(\Gamma)}b_i= \fbb^T d$. 

Let us call $f(\Gamma, x^*) = \fbb^T  \AI x^*  - t^- +  r (\onevec^T x^*  - \beta - 1)$.
It is easy to see that $x^*$ violates the $\Gamma$-inequality \eqref{eq:gamma-ineq} if and only if $f(\Gamma, x^*) <0$.
We will prove that $f(\Gamma, x^*)= c(\Gamma, x^*)$.

On one hand, we have
\begin{align}
f(\Gamma, x^*)  &= \fbb^T  \AI x^*  - t^-  +  r (\onevec^T x^*  - \beta - 1), \nonumber \\
&= \fbb^T (\AI x^* - d) - r (\beta + 1 -   \onevec^T x^*), \nonumber \\
&= \fbb^T s^* - r (\beta + 1 -   \onevec^T x^*). \label{eq:costo-gamma-3}
\end{align}
%
On the other hand, we have already seen in Lemma \ref{ppositivo} that 
$$c(\Gamma, x^*) = \fbb^T s^* + (\fb^T-\fbb^T)(\mu s^*-\mu (1-\mu)v).$$
From Corollary~\ref{th:parameter-properties} it follows that $\fbb^T s^{*} = \fb^T s^{*} - p \onevec^T x^{*} + t$. Hence, we can write
\begin{equation}
\label{eq:tau}
c(\Gamma, x^*)=\fbb^T s^*-\mu[t-(\onevec^T x^*+\mu-1)p].
\end{equation}
Similarly, if we consider that 
$$c(\Gamma, x^*) = \fb^T s^* - (\fb^T-\fbb^T) \cbb(x^*),$$
it is not hard to see that  
\begin{equation}
\label{eq:tau+uno}
c(\Gamma, x^*) = \fb^T s^*- (1-\mu)[p(\onevec^T x^*+\mu)-t].
\end{equation}

Let $\tau^*:= \floor{\onevec^T x^*}$, i.e., $\onevec^T x^* + \mu = \tau^* + 1$. Since $c(\Gamma, x^*)<0$, it follows from \eqref{eq:tau} that $t - \tau^* p > 0$ or, equivalently, $\tau^* < \frac{t}{p}$.
Similarly, from \eqref{eq:tau+uno} we obtain $t-(\tau^*+1) p < 0$, i.e., $\tau^* + 1 > \frac{t}{p}$. 
Thus, $\tau^* = \floor{\frac{t}{p}} = \beta$. 

This fact together with \eqref{eq:costo-gamma-3} and \eqref{eq:tau} imply
$$f(\Gamma, x^*) = \fbb^T s^* - (\tau^* + 1 -   \onevec^T x^*) r
= \fbb^T s^* - \mu  (t - \tau^* p) 
= c(\Gamma, x^*)<0.$$
\end{myproof}

As a consequence of the previous results we obtain a complete linear description of $\Qi{A,b}$.

\begin{theorem}
\label{th:complete-description}
For any circular matrix $A$ and any vector $b \in \Z_+^m$, the polyhedron $\Qi{A,b}$ is completely
described by the inequalities defining $\Q{A,b}$ and circuit inequalities induced by circuits $\Gamma$ in $\Aux{A}$ with $\pGam\geq 2$ such that $\pGam$ does not divide $\tGam$.
\end{theorem}

Observe that the above theorem and Lemma \ref{rem:non-redundant} imply that any facet of $\Qi{A,b}$ not in the system $Ax\geq b, x\geq 0$ must have full support. This fact has already been observed in \cite{ArgiroffoBianchi09} for the particular case of the set covering polyhedron related to circulant matrices. In contrast, there are non-boolean facets of the packing polytope of circular matrices which do not have full support. 
As a further consequence of Theorem \ref{th:complete-description} and the polynomiality of the Minimum Cost Circulation Problem, 
we also obtain that the Weighted $L$-Domination Problem is polynomial time solvable on circular interval graph. 

%
%
%

\section{The case of homogeneous right-hand side}
\label{sec:homog-rhs}

In this section, we consider polyhedra $\Qi{A,b}$ for a circular matrix $A$ and $b=\alpha \onevec$ with $\alpha \in \N$. 
Observe that for these class of polyhedra dominating rows of the matrix $A$ are associated with redundant constraints.
For this reason, in the remaining of this article we always assume that $A$ has no dominating rows.

We will prove that in this case, relevant circuit inequalities are  induced by circuits in $\Aux{A}$ without reverse row arcs.

Remind that, for every $i\in [m]$, $P_i^+$ (resp. $P_i^-$) is the path of short forward (resp. reverse) arcs connecting $l_i-1$ with $l_i+k_i-1$ (resp. $l_i+k_i-1$ with $l_i-1$).
We start with the following result.
\begin{lemma}
\label{th:path-cost}
Let $A$ be a circular matrix and 
$x^* \in \Q{A,\alpha \onevec}$ with $\alpha \in \N$. 
Then, for every $i\in [m]$ the following statements hold in $D(A,x^*)$:
\begin{enumerate}
\item[(i)]
The cost of the forward row arc $a_i=(\ell_i -1, \ell_i + k_i - 1)$ is smaller than the cost of $P^+_i$ by the amount $-\mu \alpha$; i.e., $c^+_i(x^*)-c^+(P^+_i,x^*)=-\mu \alpha$.
\item[(ii)]
The cost of the reverse row arc $\ab_i=(\ell_i + k_i - 1,\ell_i -1)$ is smaller
than the cost of $P^-_i$, by the amount $-(1 - \mu)\alpha$; i.e., $c^-_i(x^*)-c^-(P^+_i,x^*)=-(1 - \mu)\alpha$.
\end{enumerate}

\end{lemma}

\begin{myproof}
Let $i\in [m]$ and $u:= e_i - \sum_{j= \ell_i}^{\ell_i + k_i - 1} e_{m+j} \in \R^{m+n}$, where $e_k$ denotes the $k$-th canonic vector in $\R^{m+n}$.

Let us first prove that, 
$u^T s^* =u^T(\AI x^* - d)= - \alpha$ and $u^T v =0$. Indeed, 
$$
u^T \AI = u^T \binomio{A}{I}  = e_i^T A  - \sum_{j= \ell_i}^{\ell_i + k_i - 1} e_j^T  = \zerovec^T,
$$
as the $i$-th row of $A$ is the incidence vector of
$[\ell_i, \ell_i + k_i)_n \subset [n]$. Thus,
\begin{equation}
\label{us}
u^T s^* = -u^T d = -u^T  \binomio{\alpha \onevec}{0} = - \alpha.
\end{equation}

Let us now analyze the product $u^T v$. Since $v$ is the last column of the matrix $\AI$, we have that $v_i = 1$ if and only if $n \in [\ell_i, \ell_i + k_i)_n$. Furthermore, $v_{m+j}=1$ if and only if $j=n$. Hence, 
\begin{equation}
\label{uv}
u^T v = v_i - \sum_{j= \ell_i}^{\ell_i + k_i - 1} v_{m+j}= 0. 
\end{equation}
Finally, since $\cb(x^*) = \mu s^* - \mu (1- \mu) v$ and $c^+(P^+_i,x^*)= \sum_{j= \ell_i}^{\ell_i + k_i - 1} \cb_{m+j}(x^*)$, 
$$\cb_i(x^*) - c^+(P^+_i,x^*) = u^T \cb(x^*) = \mu (u^T s^*) - \mu (1 - \mu) (u^T v),$$ 
replacing (\ref{us}) and (\ref{uv}) in the last equation 
we have:  
$$\cb_i(x^*) - c^+(P^+_i,x^*) = - \mu \alpha.$$
The proof of part (ii) is similar, considering that $\cbb(x^*) = (1 - \mu) s^* + \mu (1- \mu) v$.
\end{myproof}

To prove that circuit inequalities induced by circuits with reverse row arcs are redundant, we
state at first the following result. 


\begin{lemma}
\label{th:no-reverse}
Let $A$ be a circular matrix and $\Gamma$ be a circuit in $\Aux{A}$ with $\pGam > 0$. Let $\alpha \in \N$ and $x^* \in \Q{A,\alpha \onevec}$. Then, there exists a circuit $\Gamma'$ without
reverse row arcs, such that $p(\Gamma') > 0$ and $c(\Gamma', x^*) \leq c(\Gamma, x^*)$.
\end{lemma}

\begin{myproof}

If $\Gamma$ contains no reverse row arcs then $\Gamma'=\Gamma$. 
Otherwise, let $\ab \in E_1^{-}(\Gamma)$. 

Observe that in this case $\Gamma$ must also contain at least one forward row arc, too. Indeed, from $\pGam >0$ and Lemma~\ref{th:winding-number} it follows that any node jumped by $\ab$ must also be jumped by at least two forward arcs. But then, since $\Gamma$ is simple, at least one of these have to be a row arc. 
Hence, we can choose two row arcs $\ab_i, a_r \in E(\Gamma)$, $i\neq r$,
such that $a_r$ is the first row arc preceding $\ab_i$ in $\Gamma$. Then, the circuit must
contain a simple path $P$ from $\ell_r +k_r -1$ to $\ell_i + k_i - 1$, consisting only of short arcs.
We distinguish between the two possible cases.

\emph{Case (i): $P$ contains only reverse arcs.} Figure~\ref{fig:grafo-case-i} depicts this situation.

\begin{figure}[h]
 \centering
    \includegraphics[width=0.335\textwidth]{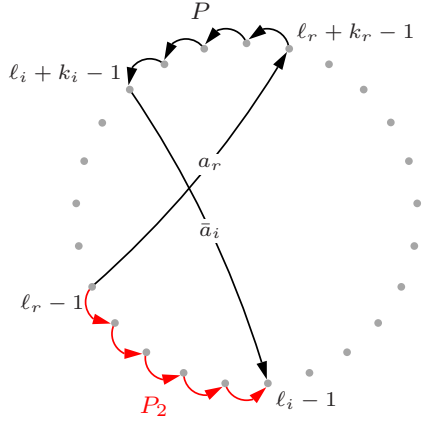}
		\caption{Case (i): $P$ consists of arcs from $E_2^{-}$.}
\label{fig:grafo-case-i}
\end{figure}

Observe that in this case $\ell_r -1$ is jumped by $\ab_i$, as otherwise, row $r$ dominates row $i$ in $A$. 
Let $P_1$ be the path from $\ell_r -1$ to $\ell_i -1$ in $\Gamma$, i.e. $P_1$ is the concatenation of 
$a_r$, $P$ and $\ab_i$. 
Consider the alternative path $P_2$ in $D(A)$ that connects $\ell_r -1$ with $\ell_i -1$  using only reverse short arcs. 

Define $\Phi$ to be the closed (not necessarily simple) path obtained from $\Gamma$ by replacing $P_1$ by $P_2$. Clearly, $\Phi$ has the same winding number and one fewer reverse row arc than $\Gamma$. Moreover,  the cost $c(\Phi, x^*)$ is smaller than or equal to the cost $c(\Gamma, x^*)$. Indeed,
\begin{align*}
c(\Gamma, x^*) - c(\Phi, x^*) &= c(P_1, x^*) - c(P_2, x^*) \\
&= \cb_r + c(P,x^*) + \cbb_i - c(P_2, x^*)\\
&= \cb_r +  \cbb_i + (c(P,x^*)- c(P_2, x^*)). \\
\end{align*}
Observe that 
$$ c(P,x^*)- c(P_2, x^*) = c(P^-_r,x^*)- c(P^-_i, x^*).$$
By Lemma \ref{th:path-cost} (ii), $c(P^-_r,x^*)- c(P^-_i, x^*)= \cbb_r-\cbb_i$. Then, 
$$
c(\Gamma, x^*) - c(\Phi, x^*)= \cb_r +  \cbb_i +  \cbb_r-\cbb_i = \cb_r+\cbb_r= s^{*}_r \geq 0. \\
$$


\emph{Case (ii): $P$ contains only forward arcs.} This situation is shown in Figure~\ref{fig:grafo-case-ii}. 

\begin{figure}[h]
 \centering
    \includegraphics[width=0.35\textwidth]{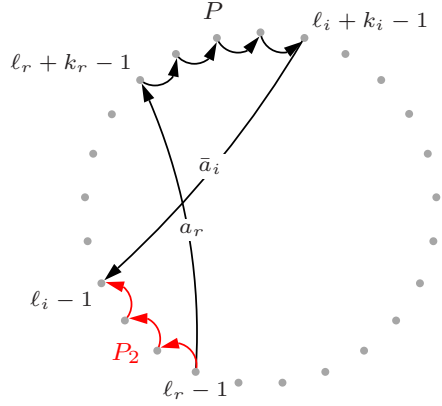}
		\caption{Case (ii): $P$ consists of arcs from $E_2^{+}$.}
\label{fig:grafo-case-ii}
\end{figure}

In this case, $\ell_i-1$ is jumped by $a_r$, as otherwise row $i$ dominates row $r$. Again, let $P_1$ be the path from $\ell_r -1$ to $\ell_i -1$ in $\Gamma$, and consider the alternative path $P_2$ in $D(A)$ consisting only from arcs of $E_2^{+}$.
Let $\Phi$ be the closed (not necessarily simple) path obtained from $\Gamma$ by replacing $P_1$ by $P_2$. This path has the same winding number, one fewer reverse row arc than $\Gamma$, and its cost $c(\Phi, x^*)$ is no larger than $c(\Gamma, x^*)$. Indeed,  
\begin{align*}
c(\Gamma, x^*) - c(\Phi, x^*) &= c(P_1, x^*) - c(P_2, x^*) \\
&= \cb_r + c(P,x^*) + \cbb_i - c(P_2, x^*)\\
&= \cb_r +  \cbb_i + (c(P,x^*)- c(P_2, x^*)). 
\end{align*}
Observe that 
$$ c(P,x^*)- c(P_2, x^*) = c(P^+_i,x^*)- c(P^+_r, x^*)$$
and from 
Lemma \ref{th:path-cost} (i) it follows that  
$c(\Gamma, x^*) - c(\Phi, x^*)=s^*_i\geq 0$.


In both cases, we have proven the existence of a closed path $\Phi$ with strictly fewer reverse row arcs than
$\Gamma$, and such that $c(\Gamma, x^*) \geq c(\Phi, x^*)$ and $p(\Gamma)=p(\Phi)>0$. But then, $\Phi$ contains at least one circuit $\Gamma^{(2)}$ with strictly fewer reverse row arcs than $\Gamma$, positive winding number, and such that $c(\Gamma, x^*) \geq c(\Gamma^{(2)}, x^*)$. 

Iterating this argument a finite number of times, we prove the existence of a circuit $\Gamma'$ without reverse
row arcs and positive winding number such that $c(\Gamma, x^*) \geq c(\Gamma', x^*)$.
\end{myproof}

As a consequence of the last lemma, 
we obtain the main result of this section:

\begin{theorem}
\label{th:result-homogeneous-case}
For any circular matrix $A$ and any $\alpha \in \N$, the polyhedron $\Qi{A,\alpha \onevec}$ is completely
described by the inequalities defining $\Q{A,\alpha \onevec}$ and circuit inequalities induced by circuits $\Gamma$ in $\Aux{A}$ without reverse row arcs, with $\pGam\geq 2$, and such that $\pGam$ does not divide $t(\Gamma, \alpha \onevec)$.

\end{theorem}

\begin{myproof}
Let $x^* \in \Q{A,\alpha \onevec}\setminus \Qi{A,\alpha \onevec}$.
Due to Lemma \ref{th:gammanonegativo} there exists at least one circuit $\Gamma$ such that $c(\Gamma, x^*) < 0$. 
From Lemma~\ref{th:no-reverse} there exists a circuit $\Gamma'$ without reverse row arcs such that $c(\Gamma', x^*) \leq c(\Gamma, x^*) < 0$. By Lemma \ref{5prima}, $x^*$ violates the $\Gamma'$-inequality. Then, 
$\Qi{A,\alpha \onevec}$ is completely
described by boolean inequalities and circuit inequalities induced by circuits $\Gamma$ in $\Aux{A}$ without reverse row arcs.
By Theorem \ref{th:complete-description}, we only need to consider $\Gamma$-inequalities such that $\pGam\geq 2$ and $\pGam$ does not divide $t(\Gamma, \alpha \onevec)$. 
\end{myproof}

In the following we further study combinatorial properties of circuits without reverse row arcs in $\Aux{A}$ and the implications for the related inequalities. In particular, in the case of the set covering polyhedron, 
we show that these circuit inequalities reduce to row family inequalities.

\section{Circuits without reverse row arcs and their inequalities}
\label{sec:no-reverse-arcs}

Given a circular matrix $A$, let us call $F(A)$ the digraph with nodes in $[n]$ and all arcs in $\Aux{A}$ except for reverse row arcs. Moreover, let $\Gamma$ be a circuit in $F(A)$. Keeping the same notation introduced in  \cite{Stauffer05}, we consider the partition of the nodes of $F(A)$ into the following three classes:
\begin{itemize}
\item[(i)] \emph{circles} 
$\oGam := \setof{j \in [n]}{(j-1,j)\in E(\Gamma)}$,
\item[(ii)] \emph{crosses} 
$\xGam :=\setof{j \in [n]}{(j, j-1)\in E(\Gamma)}$, and 
\item[(iii)] \emph{bullets} $\bGam: = [n]\setminus (\circ(\Gamma)\cup \otimes(\Gamma))$.  
\end{itemize}

Observe that circle (resp. cross) nodes are the heads (resp. tails) of forward (resp. reverse) short arcs of $\Gamma$. A bullet node is either a node outside $\Gamma$, or it is the tail or the head of a row arc.
We say that a bullet is an \emph{essential bullet} if it is reached by $\Gamma$.

Remind that a forward row arc $(u,v)$ jumps over a node $j \in [n]$ if $j \in (u  ,v]_n$. Also, the only forward (resp. reverse)  short arc jumping  over $j$ is the arc $(j-1,j)$ (resp. $(j,j-1)$).

The number of row arcs of $\Gamma$ jumping over a given node depends on which partition class it belongs to.

\begin{lemma} 
\label{th:jumping-number-nodes}
Let $A$ be a circular matrix and $\Gamma$ be a circuit in $F(A)$ with winding number $p$. For each node $j\in [n]$, let $r(j)$ be the number of row arcs of $\Gamma$ that jump over $j$.
Then, 
$$
r(j) = \left\{
\begin{array}{ll}
p-1 \; & \mbox{if $j \in \oGam$}, \\
p+1 \;& \mbox{if $j \in \xGam$}, \\
p \;& \mbox{if $j \in \bGam$}. \\
\end{array}
\right.
$$
\end{lemma}  

\begin{myproof}
From Lemma \ref{th:winding-number} we know that, for all $j \in [n]$, 
$p=p^{+}(\Gamma, j) - p^{-}(\Gamma, j)$.

If $j \in \oGam$, there is exactly one forward short arc jumping over $j$. Since $\Gamma$ is a circuit there is no reverse short arcs that jump over this node. Hence, there are exactly $p-1$ forward row arcs that jump over $j$. 

If $j \in \xGam$, again from the assumption that $\Gamma$ is a circuit there is exactly one reverse short arc and no forward short arcs that jump over $j$. It follows that $p+1$ forward row arcs must jump over this node. 

Finally, if $j \in \bGam$, neither forward nor reverse short arcs can jump over $j$ and then $r(j)= p$. 
\end{myproof}


%
From the previous results, the relevant circuit inequalities of $\Qi{A,\alpha \onevec}$ have a particular structure.

\begin{theorem}
\label{th:2-coefs}
Let $A$ be a circular matrix and $\alpha\in \N$. Let $\Gamma$ be a circuit in $F(A)$ with $s$ row arcs and winding number $p$, fulfilling the conditions of Lemma \ref{rem:non-redundant}. If 
$r:= \alpha s - p \floor{\frac{\alpha s}{p}}$, the $\Gamma$-inequality of $\Qi{A,\alpha \onevec}$ has the form:
\begin{equation}
\label{eq:2-coefs-ineq}
r \sum_{j \not\in \xGam} x_j + (r+1) \sum_{j \in \xGam} x_j \geq r\ceil{\frac{\alpha s}{p}}.
\end{equation}
Moreover, if $\alpha= 1$ and $\xGam \neq \emptyset$, this inequality is the row family inequality induced by $F:= \{i \in [m]: a_i \text{ is a row arc of } \, \Gamma\}$.
\end{theorem}

\begin{myproof}
Recall from Definition \ref{gammaineq} that the $\Gamma$-inequality has the form
$$\sum_{j\in [n]} [p^-(\Gamma,j)+r(\Gamma,b)] \, x_j \geq  r(\Gamma,b) \, (\beta (\Gamma,b) +1)+ \sum_{i\in E^-_1(\Gamma)}\!\!\! b_i.$$

Clearly, since $\Gamma$ has no reverse row arcs,  $\sum_{i\in E_1^-(\Gamma)}b_i=0$ and $t(\Gamma,b)= \sum_{i\in E_1^+(\Gamma)}b_i= \alpha s$. Then, $\beta (\Gamma, b)= \ceil{\frac{\alpha s}{p}}-1$ and $r(\Gamma, b)=r$. Then, the $\Gamma$-inequality (\ref{eq:gamma-ineq}) has the form
\begin{equation}
\sum_{j\in [n]} [p^-(\Gamma,j)+r] \, x_j \geq  r \, \ceil{\frac{\alpha s}{p}}.
\end{equation}
In order to obtain (\ref{eq:2-coefs-ineq}), it only remains to observe that, since $E^-(\Gamma)$ contains only short reverse arcs, we have: 
$$p^-(\Gamma,j)=\left\{
 \begin{array}{rl}
 1 & \mbox{ if } j \in \xGam, \\
 0 & \mbox{ otherwise.}
\end{array} \right.$$

Now assume that $\alpha= 1$ and $\xGam \neq \emptyset$. 

Since $s= \card{F}$ to prove that \eqref{eq:2-coefs-ineq} is the row family inequality induced by $F$, it suffices to show that $p= \max_{j  \in [n]} \{\sum_{i \in F} a_{ij} \} - 1$ and $O(F,p) = \xGam$. 
Indeed, it is not hard to see that, for any $j \in [n]$, $\sum_{i \in F} a_{ij}$ coincides with $r(j)$ defined in Lemma~\ref{th:jumping-number-nodes}. Moreover, since $\xGam \neq \emptyset$, 
%
$$\max_{j  \in [n]} \{\sum_{i \in F} a_{ij} \} - 1 = (p+1) - 1 =p \quad \text{ and } \quad
O(F, p) = \{j \in [n]: \sum_{i \in F} a_{ij} = p+1\} = \xGam.$$ 
Then, the $\Gamma$-inequality is the row family inequality induced by $F$ and the proof is complete.
\end{myproof}

In the particular case when $A$ is a circulant matrix, relevant circuit inequalities correspond to circuits without circle nodes.  

\begin{lemma}
Let $n,k\in \N$ such that $2\leq k\leq n-2$ and $\Gamma$ be a circuit in $F(C_n^k)$ with $s$ row arcs, winding number $p$, $\xGam\neq \emptyset$ and $\oGam\neq \emptyset$. Then, for any $\alpha\in \N$,  
the $\Gamma$-inequality is not a facet of $\Qi{C_n^k,\alpha \onevec}$. 
\end{lemma}

\begin{myproof}
Since $\xGam\neq \emptyset$ and $\oGam\neq \emptyset$, there is a path in $\Gamma$ connecting a cross with a circle. Let $P$ be a shortest path in $\Gamma$ with this condition. Assume that $P$ starts at $u\in \xGam$, has $h\geq 1$ row arcs and ends at $v\in \oGam$. Then, the nodes of $P$ are $u,u-1$, $(u-1)+jk$ with $1\leq j\leq h$, and $(u-1)+hk+1=u+hk=v$. 

Consider $P'$ the path of row arcs in $F(A)$ from $u$ to $u+hk=v$ and let $\Gamma'$ be the closed path obtained by replacing $P$ by $P'$ in $\Gamma$. We will see that $\Gamma'$ is a circuit. To do that, we only need to prove that internal nodes in $P'$ do not belong to $\Gamma$. Clearly, the internal nodes in $P'$ are $u+jk$ with $1\leq j\leq h-1$.

Assume there exists $j$ with $1\leq j\leq h-1$ such that $u+jk$ is a node of $\Gamma$. Let $t=\min \{j: u+jk\in V(\Gamma), 1\leq j\leq h-1\}$. Clearly, $u+tk\notin \oGam$ (resp. $u+tk\notin \xGam$), otherwise there are two arcs from $\Gamma$ leaving (resp. entering) $(u-1)+tk$. Then, either $(u+tk+1, u+tk)$ or $(u+(t-1)k, u+tk)$ is an arc of $\Gamma$. If $(u+tk+1, u+tk)$ is an arc of $\Gamma$, $u+tk+1\in \xGam$ and the path in $\Gamma$ from $u+tk+1$ to $v$ is shorter than $P$, a contradiction. If $(u+(t-1)k, u+tk)$ is an arc of $\Gamma$, $u+(t-1)k$ is a node of $\Gamma$. If $t=1$, we have two arcs in $\Gamma$ leaving $u$. If $t\geq 2$, we have a contradiction with the definition of $t$. 

Then, $\Gamma'$ is a circuit in $F(A)$. 
Moreover, $\Gamma'$ has $s$ row arcs, winding number $p$, and $\otimes{(\Gamma')}$ is strictly contained in $\xGam$. Hence, the $\Gamma'$-inequality implies the $\Gamma$-inequality. 
\end{myproof}

Observe that if $\Gamma$ is a circuit in $F(A)$ such that $\xGam=\emptyset$, the $\Gamma$-inequality is implied by the rank constraint.  As a consequence we have: 


\begin{theorem}
\label{circulantes}
For any positive numbers $n$ and $k$, with $2\leq k\leq n-2$ and $\alpha \in \N$, a complete linear description for the polyhedron $\Qi{C^k_n,\alpha \onevec}$ is given by
the inequalities defining $\Q{C^k_n,\alpha \onevec}$, 
the rank constraint, and 
circuit inequalities corresponding to circuits in $F(A)$ without short forward arcs. 
\end{theorem}

In the next section we will see that relevant inequalities for the set covering polyhedron of circular matrices are minor induced row family inequalities.

\section{Set covering polyhedron of circular matrices and circulant minors}
\label{sec:set-covering-polyhedron}

Throughout this section, we restrict our attention to the set covering polyhedron of circular matrices. Remind that we have assumed that the matrix $A$ has no dominating rows.

Let $\Gamma$ be a circuit of $F(A)$ with winding number $p$ and $s$ essential bullets $\{b_j: j=1,\ldots,s\}$, with $1\leq b_1 < b_2<\ldots< b_s\leq n$. 

Clearly, if $[u,w]_n\subset \oGam$  and  $u-1,w+1\notin \oGam$ (resp. $[u,w]_n\subset \xGam$ and $u-1,w+1\notin \xGam$) then $u-1$ is an essential bullet and $w+1\in \bGam$.

\begin{figure}[h]
 \centering
  \begin{minipage}[c]{4cm}
   \label{f:circleblock}
    \includegraphics[width=0.7975\textwidth]{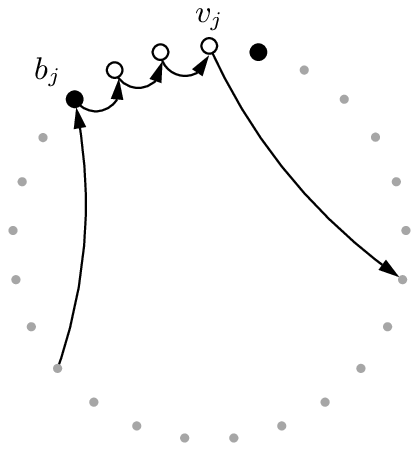}
    \end{minipage}
  \begin{minipage}[c]{4cm}
   \label{f:crossblock}
    \includegraphics[width=0.7975\textwidth]{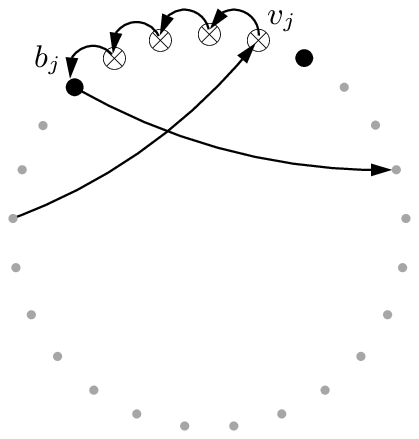}
    \end{minipage}
  \begin{minipage}[c]{4cm}
   \label{f:bulletblock}
    \includegraphics[width=0.7975\textwidth]{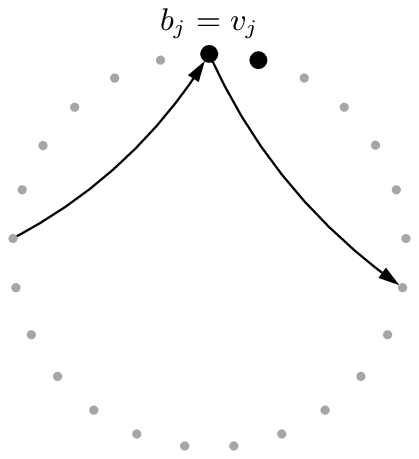}
    \end{minipage}
 \caption{All three possible block types: (a) circle block, (b) cross block, (c) bullet block.}
 \label{fig:block-types}
\end{figure}

For each $j\in [s]$ we define $v_j$ as the node of $\Gamma$ in $[b_j,b_{j+1})$ such that:
\begin{itemize}
\item[(i)] if $b_j+1 \in \oGam$ then $[b_j+1,v_j]_n\subset \oGam$ and $v_{j}+1\notin \oGam$, 
\item[(ii)] if $b_j+1 \in \xGam$  then $[b_j+1,v_j]_n\subset \xGam$ and $v_{j}+1\notin \xGam$,
\item[(iii)] if $b_j+1 \in \bGam$  then $v_j=b_j$.
\end{itemize}
Then, for each $j\in [s]$ we define the \emph{block} $B_j=[b_j, v_j]_n$ which can be a \emph{circle} block, a \emph{cross} block or a \emph{bullet} block depending on if $b_j+1$ is a circle, a cross or a bullet of $\Gamma$. 
It is easy to check that the blocks $\{B_j: j=1,\ldots,s\}$ define a partition of nodes of $\Gamma$.
Figure~\ref{fig:block-types} illustrates these three type of blocks.

\begin{remark}\label{arcos} 
For each $j\in [s]$ there exists one row arc leaving $B_j$ and another row arc entering $B_j$. 
Let $B_j^{-} \in B_j$ be the tail of the arc leaving $B_j$, while $B_j^{+} \in B_j$ denotes the head of the arc entering in $B_j$. In particular, if $B_j$ is a cross block, $B^-_j=b_j$ and $B^+_j=v_j$. If $B_j$ is a circle block, $B^-_j=v_j$ and $B^+_j=b_j$. Finally, if $B_j$ is a bullet block, $B^-_j=B^+_j=b_j=v_j$.

Observe that if $B_j$ is a circle block (cross block) then there is always a path of forward (reverse) short arcs in $\Gamma$ that joins $B_j^{+}$ with $B_j^{-}$. Moreover, the nodes of every path in $\Gamma$ consisting only of short arcs belong to the same block. 
\end{remark}

We have the following result:

\begin{lemma}
\label{peb}
Let $A$ be a circular matrix, $\Gamma$ be a circuit of $F(A)$ with winding number $p$ and $s$ essential bullets $\{b_j: j\in [s]\}$, with $1\leq b_1 < b_2<\ldots< b_s\leq n$. Then, $\Gamma$ has $s$ row arcs and each row arc in $\Gamma$ jumps over $p$ essential bullets, i.e., it has the form $(B^-_i,B^+_{i+p})$ for some $i\in [s]$. Moreover, $\gcd(s,p)=1$.
\end{lemma}
\begin{myproof}
Since there are $s$ blocks and each row arc of $\Gamma$ is the leaving arc of exactly one block, $\Gamma$ has $s$ row arcs. 

Consider the row arc $(B^-_i,B^+_{i+t})$ of $\Gamma$.  Clearly, it jumps over the $t$ essential bullets $b_{i+\ell}$ with $1\leq \ell \leq t$.
Moreover, since $A$ has no dominating rows, 
then $(B^-_{i+1},B^+_{i+1+t})$ is also a row arc of $\Gamma$. Iterating this argument, one can verify that the row arcs of $\Gamma$ jumping over $b_{i+t}$ are exactly $(B^-_{i+\ell},B^+_{i+\ell+t})$ with $0\leq \ell \leq t-1$. From Lemma \ref{th:jumping-number-nodes}, it follows that $t=p$.
Thus, $(B^-_i,B^+_{i+t})$ jumps over $p$ essential bullets. 

%


Let $\BGra=(\BVer, \BEdg)$ be the directed graph where $\BVer=[s]$ and $(i,j)\in \BEdg$ if $(B^-_i, B^+_j)$ is a row arc of $\Gamma$. Hence, $j=i+p$ and $\BGra$ is a circuit 
with $s$ arcs of length $p$. 
But then, 
$\gcd(s,p)=1$ must holds. 
\end{myproof}

Lemma \ref{peb} also establishes a relationship between circuits in $F(A)$ and some circulant submatrices of $A$. 

\begin{corollary}
\label{circulantsubmatrix}
Let $A$ be a circular matrix and $\Gamma$ be a circuit in $F(A)$ with winding number $p$ and $s$ row arcs. Let $L\subset [n]$ be the set of essential bullets of $\Gamma$ and
$F\subset [m]$ be the set of rows of $A$ corresponding to the row arcs of $\Gamma$.
Then, the submatrix of $A$ induced by rows in $F$ and columns in $L$ is isomorphic to the circulant matrix $\C{s}{p}$ with $\gcd(s,p)=1$. 
\end{corollary}

\begin{myproof}
Let $A'$ be the submatrix of $A$ induced by the rows in $F$ and the columns in $L$.
From Lemma \ref{peb}, all the row arcs in $\Gamma$ jump over $p$ essential bullets. Then, each row of $A'$ has exactly $p$ entries equal to one and there is no pair of equal rows in $A'$. Hence, $A'$ is isomorphic to $\C{s}{p}$.
\end{myproof}

Observe that the submatrix mentioned in the corollary above is not necessarily a circulant (contraction) minor of 
$A$. Indeed, after deleting the columns from $N=[n]\setminus L$, there might be rows in $[m]\setminus F$ that
are dominated by rows in $F$. This only happens when there is a row arc in $F(A)$ that jumps over less than $p$ essential bullets of $\Gamma$. Inspired by the results in \cite{Stauffer05}, we say that a row arc in $F(A)$ is a \emph{bad arc (with respect to $\Gamma$)} if it jumps over less than $p$ essential bullets of $\Gamma$.

Then, we have:

\begin{corollary}
\label{th:circulant-minor}
Let $A$ be a circular matrix and $\Gamma$ be a circuit in $F(A)$ with winding number $p$ and $s$ row arcs. Let $L\subset [n]$ be the set of essential bullets of $\Gamma$ and $F\subset [m]$ be the set of rows of $A$ corresponding to the row arcs of $\Gamma$. 
If $\Gamma$ has no bad arcs, then the minor $A/N$ of $A$ is isomorphic to the circulant matrix $\C{s}{p}$. Moreover, the $\Gamma$-inequality for $\Qi{A}$ is a minor related row family inequality. 
\end{corollary}

The following result gives a characterization of bad arcs in terms of their endpoints:

\begin{theorem}
\label{th:arcs-vs-ess-bullets}
Let $A$ be a circular matrix and $\Gamma$ be a circuit of $F(A)$ with $s$ essential bullets $1\leq b_1<\ldots<b_s\leq n$, and winding number $p$. Let $(u,v)$ be a row arc in $F(A)$ that jumps over $k$ essential bullets of $\Gamma$. Then, $k \in \{p-1, p, p+1\}$. Moreover, $(u,v)$ jumps over $p-1$ essential bullets of $\Gamma$ if and only if the following two conditions hold:
\begin{itemize}
	\item[(i)] $u$ belongs to a circle block of $\Gamma$ 
	\item[(ii)] $v\in \oGam$ or $v$ is not a node of $\Gamma$.
\end{itemize}
In addition, if $u\in B_i$ then $u\neq v_i$. In this case, if $v\in \oGam$, $B_{i+p-1}$ is a circle block and $v\in B_{i+p-1}\setminus \{b_{i+p-1}\}$. If $v$ is not a node of $\Gamma$, $v\in (v_{i+p-1},b_{i+p})_n$.
\end{theorem} 

\begin{myproof}
Due to Lemma \ref{peb} if $(u,v)$ is a row arc of $\Gamma$, $k=p$. Now consider a row arc $(u,v)$ not in $\Gamma$. 


If $u$ is not a node of $\Gamma$, then it is between two consecutive blocks in $\Gamma$, i.e., there exists $i\in [s]$ such that $u\in (v_{i-1}, b_i)_n$. Since there are no dominating rows in $A$ and $(B^-_{i-1}, B^+_{i+p-1})$,  $(B^-_{i}, B^+_{i+p})$ are row arcs of $\Gamma$ then $v\in (b_{i+p-1},v_{i+p})_n$. But then, $(u,v)$ jumps either over $p$ essential bullets when $v\in (b_{i+p-1},b_{i+p})_n$, or over $p+1$ essential bullets when $v\in [b_{i+p},v_{i+p})_n$. 
Hence, if $(u,v)$ is a bad arc, $u$ has to be a node of $\Gamma$.

Now assume $u$ is a node of $\Gamma$ and $u\in B_i$, for some $i\in [s]$. Again, since there is no dominating rows in $A$, and $(B^-_{i-1}, B^+_{i+p-1})$, $(B^-_{i+1}, B^+_{i+p+1})$ are row arcs of $\Gamma$ then $v\in (b_{i+p-1},v_{i+p+1})_n$. Therefore, $(u,v)$ jumps either over $p-1$, $p$ or $p+1$ essential bullets, depending on whether $v\in (b_{i+p-1},b_{i+p})_n$, $v\in [b_{i+p},b_{i+p+1})_n$ or $v\in [b_{i+p+1},v_{i+p+1})_n$, respectively.

Thus, any row arc in $F(A)$ jumps over $k$ essential bullets of $\Gamma$ with $k\in \{p-1,p,p+1\}$. Moreover, $(u,v)$ is a bad arc if it jumps over exactly $p-1$ essential bullets. In this case, $u$ is a node of $\Gamma$ and if $u\in B_i$, $v\in (b_{i+p-1},b_{i+p})_n$. Let us analyze this last case.

Since $v\in (b_{i+p-1},b_{i+p})_n$ then $B^+_i=v_i$, as otherwise the row of $A$ corresponding to the arc $(B^-_i,B^+_{i+p})$ is dominated by the row corresponding to $(u,v)$. Therefore, $B_i$ is a circle block and $u\neq v_i$. Moreover, if $v$ is a node of $\Gamma$ then $B^+_{i+p-1}=b_{i+p-1}$, as otherwise the row corresponding to the arc  $(B^-_{i-1},B^+_{i+p-1})$ is dominated by the row corresponding to $(u,v)$. As a consequence, $B_{i+p-1}$ is a circle block and $v\in \oGam$.
\end{myproof}

\begin{figure}[h]
 \centering
   \begin{minipage}[c]{2cm}
    \end{minipage}
  \begin{minipage}[c]{6cm}
   \includegraphics[width=.78\textwidth]{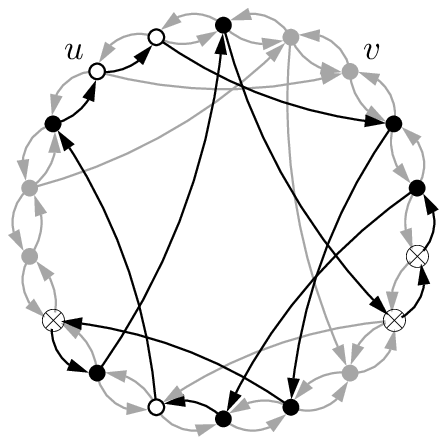}
    \end{minipage}
  \begin{minipage}[c]{6cm}
    \includegraphics[width=.78\textwidth]{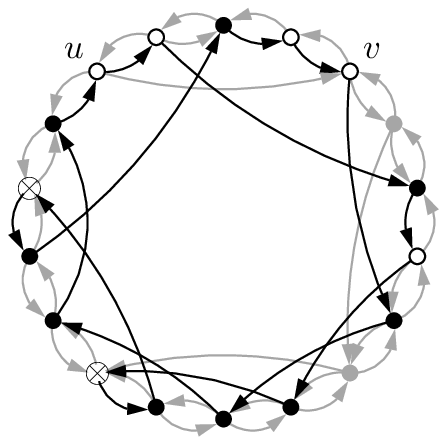}
    \end{minipage}
 \caption{(a) Case $v$ is not a node of $\Gamma$. (b) Case $v$ is a circle of $\Gamma$.}
 \label{fig:badarc}
\end{figure}

Figure~\ref{fig:badarc} depicts the two possible situations for a row arc that jumps over one essential bullet in a circuit with winding number two.  

%

Clearly, if $\oGam=\emptyset$, there is no bad arc with respect to $\Gamma$. Then, as a consequence of Theorem \ref{circulantes} and Corollary \ref{th:circulant-minor} we obtain the following result that was conjectured in \cite{Torres15}:

\begin{corollary}
For any positive numbers $n$ and $k$, with $2\leq k\leq n-2$ a complete linear description for the set covering polyhedron $\Qi{C^k_n}$ is given by boolean inequalities, the rank constraint, and minor related row family inequalities induced by circulant minors $\C{s}{p}$ of $\C{n}{k}$, with $\gcd(s,p)=1$.
\end{corollary}

In addition, recall that any circulant matrix $C_n^k$ with $k$ odd corresponds to the closed neighborhood matrix of a web graph and conversely. Then, the last result yields a complete description of the dominating set polyhedron of web graphs by boolean inequalities and row family inequalities induced by circulant minors $C_s^p$ with $\gcd(s,p)=1$. 
Moreover, the corollary above can be seen as the counterpart of the complete description of the stable set polytope of web graphs given in \cite{Stauffer05}, proving a previous conjecture stated in \cite{PecherWagler06}.
Therein, the polytope is obtained by clique inequalities and clique family inequalities associated with subwebs $W_s^{p-1}$ with $\gcd(s,p)=1$. 

\medskip

In the remaining of this section, we will see that minor related row family inequalities are sufficient for describing the set covering polyhedron of any circular matrix. In order to do so, we prove that for every circular matrix $A$, the inequalities induced by circuits without bad arcs are sufficient for describing the set covering polyhedron $\Qi{A}$. 
 

Let $x$ and $y$ be two nodes of $\Gamma$ that belong to a same block $B_i$, for some $i\in [s]$. Denote by $\Pi (x,y)$ the path of short arcs in $\Gamma$ that goes from node $x$ to node $y$. If $x=y$ then $\Pi (x,y)$ is the emptyset. Observe that $\Pi(x,y)$ is contained in $\Pi(B_i^{+},B_i^{-})$. 

In addition, if $x$ and $y$ are two distinct nodes in $[n]$, let $\pi (x,y)$ be the path of short forward arcs in $F(A)$ that goes from node $x$ to node $y$. It is clear $\pi (x,y)$ is nonempty and simple.

Let $(u,v)$ be a bad arc with respect to $\Gamma$. From Theorem \ref{th:arcs-vs-ess-bullets} it holds that $u$ belongs to a circle block and $v$ is either a node in another circle block or it is outside $\Gamma$. 

We will first see that if $v$ is a node of $\Gamma$, the $\Gamma$-inequality is not a facet of $\Qi{A}$ or it coincides with a $\Gamma'$-inequality where $\Gamma'$ is a circuit in $F(A)$ having less bad arcs than $\Gamma$. 
For this purpose, we need some technical previous results. 

Assume w.l.o.g. that $u\in B_1$ and consequently $B_1$ is a circle block. Since $v$ is a node of $\Gamma$, from Theorem \ref{th:arcs-vs-ess-bullets}, $v \in B_p$ and $B_p$ is a circle block. Moreover, $u \in[b_1,v_1)_n$ and $v\in(b_p,v_p]_n$.

\begin{figure}[h]
 \centering
   \begin{minipage}[c]{2cm}
    \end{minipage}
  \begin{minipage}[c]{5.5cm}
   \includegraphics[width=0.825\textwidth]{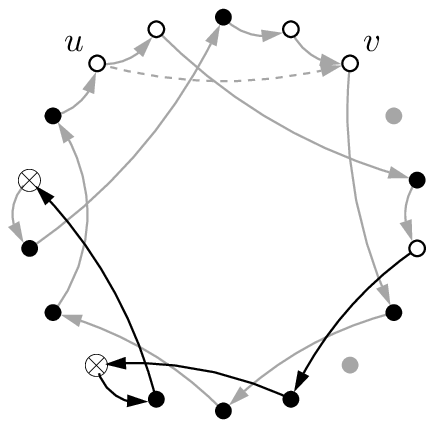}
    \end{minipage}
  \begin{minipage}[c]{5.5cm}
    \includegraphics[width=0.825\textwidth]{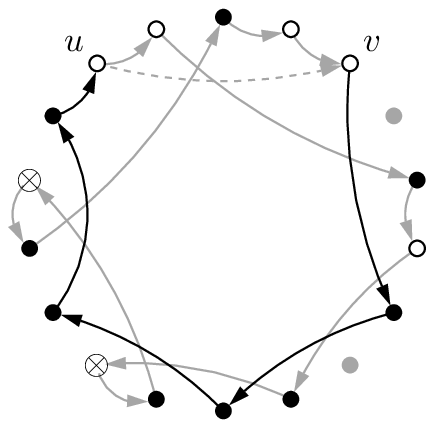}
    \end{minipage}
 \caption{Paths $P_1$ (left figure, black lines), $P_2$ (rigth figure, black lines) in a circuit $\Gamma$ (continuous lines) with a bad arc $(u,v)$ (dashed line).}
 \label{fig:caminosbadarcb}
\end{figure}

Let $P_1$ be the  path in $\Gamma$ that goes from $B_{p+1}^{-}$ to $B_s^{+}$ and $P_2$ be the path in $\Gamma$ that goes from $v$ to $u$. An example is illustrated in Figure \ref{fig:caminosbadarcb}. 

Clearly, $\Gamma$ can be seen as the concatenation of 
$P_2$, 
$\Pi(u, v_1)$, the row arc $(v_1, B^{+}_{p+1})$, 
$\Pi(B_{p+1}^{+},B_{p+1}^{-})$, 
$P_1$, $\Pi (B^{+}_{s}, B^{-}_{s})$,  
the row arc $(B_s^-,b_p)$ and 
$\Pi (b_p,v)$. It follows that $P_1$ and $P_2$ are node disjoint paths. 

\begin{figure}[h]
 \centering
    \begin{minipage}[c]{2cm}
    \end{minipage}
  \begin{minipage}[c]{5.5cm}
   \includegraphics[width=0.825\textwidth]{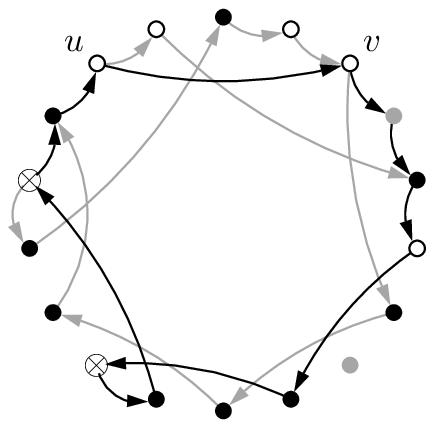}
    \end{minipage}
  \begin{minipage}[c]{5.5cm}
    \includegraphics[width=0.825\textwidth]{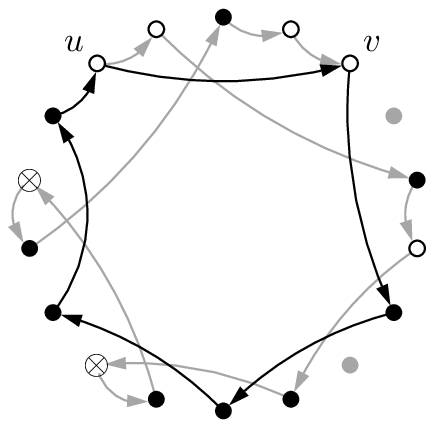}
    \end{minipage}
\caption{Circuits $\Gamma_1$ (left) and $\Gamma_2$ (right).}
 \label{fig:circuitosbadarcb}
\end{figure}

Define $\Gamma_1$ as the circuit in $F(A)$ obtained by joining $P_1$ with 
$\pi(B_{s}^{+},u)$, the row arc $(u, v)$, and 
$\pi(v,B_{p+1}^{-})$. Similarly, we define $\Gamma_2$ as the circuit obtained by joining $P_2$ together with the row arc $(u, v)$ (see Figure \ref{fig:circuitosbadarcb}).

Observe that, in $\Gamma_1$, nodes $v$ and $B_{s}^{+}$ are essential bullets but $u$ and $B_{p+1}^{-}$ are circles.  Meanwhile, the internal nodes of path $P_1$ belong to the same class with respect to $\Gamma_1$ and to $\Gamma$.

Concerning $\Gamma_2$, $v$ is an essential bullet and $u$ together with the internal nodes of $P_2$ belong to the same class with respect to $\Gamma_2$ and to $\Gamma$. Hence, every node in the set $\tabulatedset{b_2, \ldots, b_{p-1}}$ is an essential bullet for either $\Gamma_1$ or $\Gamma_2$.

In the following, it is convenient to make the next assumptions:
 
\begin{assumption}\label{ass}
$A$ is a circular matrix, $\Gamma$ is a circuit in $F(A)$ with $s$ row arcs and winding number $p\geq 2$, $\beta=\left\lfloor \frac{s}{p}\right\rfloor$ and $r=s-p\beta\geq 1$. Moreover, $(u,v)$ is a bad arc with respect to $\Gamma$ such that $u$ and $v$ belong to the circle blocks $B_1$ and $B_p$ of $\Gamma$, respectively. For $i=1,2$, $P_i$ and $\Gamma_i$ are the path and the circuit, respectively, defined above. Finally, $s_i$ and $p_i$ denote, respectively, the number of row arcs and the winding number of $\Gamma_i$, $\beta_i=\floor{\frac{s_i}{p_i}}$, and $r_i=s_i-\beta_i p_i$. 
\end{assumption}

Circuits $\Gamma_1$ and $\Gamma_2$ satisfy the following properties:

\begin{lemma} 
Under Assumption \ref{ass}, it holds:
\label{rem:crosses}
\begin{itemize}
	\item[(i)] $\xGamu \cap \xGamd = \emptyset$ and $\xGamu \cup \xGamd \subseteq \xGam$.     
	\item[(ii)] $s= s_1 + s_2$.   
\item[(iii)] $p=p_1+p_2$.   
\end{itemize}
\end{lemma}

\begin{myproof}
Let $i=1,2$. Then:
\begin{itemize}
	\item[(i)]  Since $P_1$ and $P_2$ are disjoint, $\xGamu \cap \xGamd = \emptyset$. Moreover, every cross node of $\Gamma_i$ is an internal node of $P_i$ and then it is a cross node of $\Gamma$.      
	\item[(ii)] Clearly every row arc of $\Gamma_i$ different from $(u, v)$ is a row arc of $\Gamma$ contained in the path $P_i$. Besides, $\Gamma$ contains the two row arcs $(B_s^{-},b_p)$ and $(v_1,B_{p+1}^{+})$, which are neither in $\Gamma_1$ nor in $\Gamma_2$.  
\item[(iii)] From construction, $(u,v)$ jumps over each node in the set $\tabulatedset{b_2, \ldots, b_{p-1}}\cup \{v\}$. Let $R_i= \bGami \cap \tabulatedset{b_2, \ldots, b_{p-1}}$. Since $(u, v)$ belongs to $\Gamma_i$, by  Lemma~\ref{peb} $(u, v)$ jumps over $p_i$ bullets of $\Gamma_i$: one of them is $v$ and the other $p_i-1$ bullets belong to $\tabulatedset{b_2, \ldots, b_{p-1}}$. Thus, $p_i= \card{R_i} + 1$.

It is clear that $R_1 \cap R_2 = \emptyset$  and $R_1 \cup R_2 =\tabulatedset{b_2, \ldots, b_{p-1}}$. 
Then,  $p_1+p_2=\card{R_1}+\card{R_2}+ 2= |\tabulatedset{b_2, \ldots, b_{p-1}}|+2=p$.  
\end{itemize}
\end{myproof}



\begin{lemma}
\label{th:gamma-impliedby-gam1-gam2}
Under Assumption \ref{ass}, it holds that:
\begin{itemize}
\item[(i)] if $r=1$ then 
the $\Gamma$-inequality is not facet defining for $\Qi{A}$, or it coincides with the $\Gamma_2$-inequality, 
\item[(ii)] if $r=p-1$ then the $\Gamma$-inequality is not facet defining for $\Qi{A}$.
\end{itemize} 
\end{lemma}

\begin{myproof}
\begin{itemize}
\item [(i)] Since $r=1$, the path $P_1$ connecting $B_{p+1}^{-}$ with  $B_s^{+}$ in $\Gamma$ contains exactly $\beta -1$ row arcs. Thus, from construction, $s_1=\beta$ and $p_1=1$. 
Therefore, by Lemma \ref{rem:crosses}, the circuit $\Gamma_2$ contains exactly $s_2=s - s_1 =s-\beta=\beta(p-1) + 1$ row arcs and $p_2= p-1$. 

Finally, we have $\ceil{\frac{s_2}{p_2}}= \beta + 1 = \ceil{\frac{s}{p}}$, $r_2= s_2 - p_2 \floor{\frac{s_2}{p_2}} = 1 = r$, and $\xGamd \subseteq \xGam$. As a consequence, if $\xGamd \subset \xGam$ then the $\Gamma$-inequality cannot be facet defining, as it is implied by the stronger $\Gamma_2$-inequality. Otherwise, both inequalities coincide.
\item [(ii)]
Since $r=p-1$, the path $P_2$ connecting $v$ with $u$ in $\Gamma$  
contains exactly $\beta$ row arcs and hence $s_2=\beta+1$ and $p_2=1$.
Then, by Lemma \ref{rem:crosses}, $p_1= p-1$ and $s_1=(p-1)(\beta+1)-1$, implying $\ceil{\frac{s_1}{p_1}}= \beta + 1 = \ceil{\frac{s}{p}}$ and $r_1= s_1 - p_1 \floor{\frac{s_1}{p_1}} = p-2 = r-1$. 

The $\Gamma_1$-inequality has the form:
\begin{equation}
\label{eq:gamma1-ineq}
(r-1) \!\! \sum_{j \not\in \xGamu} x_j + r \!\! \sum_{j \in \xGamu} x_j \geq (r-1)\ceil{\frac{s}{p}}.
\end{equation}
On the other hand, if we add the $s_2$ inequalities from $\Q{A}$ corresponding to the row arcs of $\Gamma_2$,
by Lemma~\ref{th:jumping-number-nodes}, we obtain the following valid inequality for $\Qi{A}$:
\begin{align*}
(p_2-1) \!\! \sum_{j \in \oGamd} x_j + p_2 \!\! \sum_{j \in \bGamd} x_j  + (p_2+1) \!\! \sum_{j \in \xGamd} x_j  &\geq s_2 = \beta + 1 = \ceil{\frac{s}{p}}
\end{align*}
which implies the valid inequality:
\begin{align*}
p_2 \!\! \sum_{j \not\in \xGamd} x_j + (p_2 + 1)\!\! \sum_{j \in \xGamd} x_j &\geq \ceil{\frac{s}{p}}.
\end{align*}

Since $p_2=1$ and, from Lemma \ref{rem:crosses} (i), we have $\xGamu \cap \xGamd = \emptyset$ and $\xGamu \cup \xGamd \subseteq \xGam$, it follows that the $\Gamma$-inequality is implied by the sum of \eqref{eq:gamma1-ineq} and the last inequality. Hence, the $\Gamma$-inequality is not facet defining. 
\end{itemize}
\end{myproof}

It remains to consider the case $r\in \tabulatedset{2, \ldots, p-2}$.

\begin{lemma}
\label{th:gamma-sum-gam1-gam2}
If Assumption \ref{ass} and $2\leq r\leq p-2$ hold then
$\beta_1=\beta_2=\beta$ and  
$r_1 + r_2 = r$. 
\end{lemma}

\begin{myproof}
As in the proof of Lemma \ref{rem:crosses}(iii), for $i=1,2$, let $R_i= \bGami \cap \tabulatedset{b_2, \ldots, b_{p-1}}$. Recall that $p_i= \card{R_i} + 1$. Since $2\leq r \leq p-2$ then $b_r, b_{r+1} \in R_1\cup R_2$. Define $R_i^{-}:= \setof{b_\ell \in R_i}{\ell < r}$ and $R_i^{+}:= \setof{b_\ell \in R_i}{\ell > r+1}$. 

Observe that, from Lemma \ref{peb}, for every $\ell \in [s]$ and every $\alpha \in \Z_+$, the path in $\Gamma$ starting at $B_{\ell}^{-}$ that uses $\alpha$ row arcs, arrives at the block $B_{\ell+\alpha p}$.
It follows that the path in $\Gamma$ starting at node $B_\ell^{-}$ for some $\ell \in  R_i^{-}$ reaches the node $B_{\ell+p-r}^{+}$ after $\beta + 1$ row arcs. It is straightforward to verify that $\ell+p-r \in \tabulatedset{2, \ldots, p-1}$. In addition, a path in $\Gamma$ starting at node $B_\ell^{-}$ for some $\ell \in  R_i^{+}$ reaches the node $B_{\ell-r}^{+}$ after $\beta$ row arcs. Since $r+2 \leq \ell \leq p-1$ and $r \geq 2$ it follows that $2 \leq \ell - r \leq p-3$. In both cases, such  paths may belong to either $\Gamma_1$ or $\Gamma_2$.

Now, the path in $\Gamma$ starting at the node $B_r^{-}$ reaches the node $B_s^{+}$ after $\beta$ row arcs, since $\beta p+r =s$. Hence, this path belongs to $\Gamma_1$. 
It follows that $b_{r}\in R_1$, $\card{R_1} \geq 1$, and $\card{R_1^{-}}< \card{R_1}$.
Moreover, from construction, $\Gamma_1$ continues with the arc $(u,v)$ and then, by following a path of short arcs, it reaches $B_{p+1}^{-}$. From this last node and after $\beta$ row arcs it reaches $B_{p+1-r}^+$. It is clear that $3\leq p+1-r\leq p-1$. Hence, the path in $\Gamma_1$ that connects $B^-_r$ with  $B_{p+1-r}^+$ contains $2\beta +1$ row arcs. 
   
Finally, the path in $\Gamma$ that starts at the node $B_{r+1}^{-}$ reaches the node $B_1^{+}$ after $\beta$ row arcs. Since $B_1$ is a circle block, $B_1^{+}=b_1$ and thus the path belongs to $\Gamma_2$. Thus, $b_{r+1}\in R_2$, $\card{R_2} \geq 1$, and $\card{R_2^{-}}< \card{R_2}$. Moreover, the circuit $\Gamma_2$ continues with $(u,v)$ and another path of short arcs until it reaches $B_p^{-}=v_p$. From this node and after $\beta$ row arcs, $\Gamma_2$ reaches $B_{p-r}^+$. Hence, the path that connects $B_{r+1}^{-}$ with $B_{p-r}^+$ in $\Gamma_2$ contains $2\beta +1$ row arcs. 

For $i=1,2$, let $\mathcal{P}_i$ be the set of simple directed paths obtained by splitting $\Gamma_i$ at the nodes
in $R_i$, i.e., the end nodes of each path in $\mathcal{P}_i$ belong to $R_i$ and no node of $R_i$ is an internal node of the path. As we have just observed, $\card{R_i} \geq 1$ holds.
If $\card{R_i} =1$, $\mathcal{P}_i$ contains one closed path, which coincides with $\Gamma_i$.  
Hence, the number of row arcs of $\Gamma_i$ can be computed by adding up the number of row arcs in each path in $\mathcal P_i$:
\begin{align*}
s_i &= \card{R_i^{-}} (\beta + 1) + \card{R_i^{+}} \beta + 2\beta + 1 \\
&= (\card{R_i^{-}} + \card{R_i^{+}} + 2) \beta + \card{R_i^{-}} + 1 \\
&= p_i \beta + \card{R_i^{-}} + 1.
\end{align*}

Observe that $1 \leq \card{R_i^{-}} + 1 \leq p_i - 1$. Hence, $\floor{\frac{s_i}{p_i}}= \beta$ and $r_i = s_i - p_i \beta= \card{R_i^{-}} + 1$. But then, $r_1 + r_2 =  \card{R_1^{-}} + \card{R_2^{-}} + 2 = \card{\tabulatedset{b_2, \ldots, b_{r-1}}} + 2 = r$.  Similarly, $\beta_i= \beta = \floor{\frac{s}{p}}$. 
\end{myproof}

As a consequence of the previous lemma we have: 

\begin{corollary}
\label{cor:gamma-sum-gam1-gam2}
Under Assumption \ref{ass} and $2\leq r\leq p-2$ hold, the $\Gamma$-inequality is not facet defining for $\Qi{A}$. 
\end{corollary}

\begin{myproof}
For $i \in \tabulatedset{1, 2}$, the $\Gamma_i$-inequality has the form
$$r_i \sum_{j \not\in \xGami} x_j + (r_i +1) \sum_{j \in \xGami} x_j \geq r_i \ceil{\frac{s_i}{p_i}}.$$
Adding the inequalities corresponding to $\Gamma_1$ and $\Gamma_2$, from Lemma \ref{th:gamma-sum-gam1-gam2}
together with Lemma~\ref{rem:crosses}(i), we have:
\begin{align*}
r \sum_{j \not\in \xGam} x_j + (r +1) \sum_{j \in \xGam} x_j 
&\geq
(r_1 + r_2) \!\!\!\!\!\!\!\! \sum_{j \not\in \xGamu \cup \xGamd} \!\!\!\!\!\!\!\! x_j + (r_1 + r_2 +1) \!\!\!\!\!\!\!\! \sum_{j \in \xGamu \cup \xGamd} \!\!\!\!\!\!\!\! x_j  \\
&=
\sum_{i=1}^2 \left( r_i \sum_{j \not\in \xGami} x_j + (r_i +1) \sum_{j \in \xGami} x_j \right) \\
&\geq
\sum_{i=1}^2  r_i \ceil{\frac{s_i}{p_i}} = (r_1 + r_2) \ceil{\frac{s}{p}} = r \ceil{\frac{s}{p}}.
\end{align*}
\end{myproof}


Finally we obtain the following result: 
\begin{theorem}
\label{th:v-in-V-redundant}

Under Assumption \ref{ass}, if  
the $\Gamma$-inequality is 
facet defining for $\Qi{A}$, then there exists a circuit 
in $F(A)$ defining the same circuit inequality and having less bad arcs than $\Gamma$.
\end{theorem}

\begin{myproof}
Assume that the $\Gamma$-inequality is facet defining for $\Qi{A}$.
By Lemma \ref{th:gamma-impliedby-gam1-gam2}(ii) and Corollary~\ref{cor:gamma-sum-gam1-gam2},  
$r=1$. 
By Lemma \ref{th:gamma-impliedby-gam1-gam2}(i), the $\Gamma$-inequality 
coincides with the $\Gamma_2$-inequality. We will prove that $\Gamma_2$ has less bad arcs than $\Gamma$. Clearly, $(u,v)$ is a bad arc for $\Gamma$ which is not a bad arc for $\Gamma_2$.  
 Thus, it suffices 
to prove that every  bad arc for $\Gamma_2$ is a bad arc for $\Gamma$. 

Assume there is a bad arc $(u',v')$ for $\Gamma_2$ which is not a bad arc for $\Gamma$. Since $\Gamma_2$ has winding number $p-1$, $(u',v')$ jumps over $p-2$ essential bullets of $\Gamma_2$ and at least $p$ essential bullets of $\Gamma$. Then, $(u',v')$ must jump over at least two essential bullets of $\Gamma$ that are not essential bullets of $\Gamma_2$. By construction of $\Gamma_2$, the only essential bullets of $\Gamma$ that are not essential bullets of $\Gamma_2$ are the nodes in the set $S=\{b_p\}\cup \{b_{1+tp}: t=1,\ldots, \beta-1\}$. The only pair in $S$ that can be jumped over by the same row arc is the pair $b_p, b_{p+1}$. But, if $(u',v')$ jumps over this pair of nodes, it must also jump over $v$, as $v\in (b_p,b_{p+1})_n$. Finally, since $v$ is an essential bullet of $\Gamma_2$, but not an essential bullet of $\Gamma$, $(u',v')$ must jump over a third essential bullet in $S$, which is not possible. 
\end{myproof}

\medskip

Now consider the case where $(u,v)$ is a bad arc with respect to $\Gamma$
and $v\notin V(\Gamma)$. 
The following result holds.

\begin{lemma}
\label{th:v-not-in-VV}
Let $A$ be a circular matrix and $\Gamma$ be a circuit in $F(A)$ with winding number $p\geq 2$, $s$ row arcs and essential bullets 
$1\leq b_1<\ldots \leq b_s\leq n$. Let $(u,v)$ be a bad arc with respect to $\Gamma$ with $u\in B_i=[b_i,v_i)_n$ and 
$v\in (v_{i+p-1}, b_{i+p})_n$, for some $i\in [s]$. If the block $B_{i+p}$ is a cross block, the $\Gamma$-inequality is not a facet defining inequality for $\Qi{A}$. 
\end{lemma}

\begin{myproof}
Assume w.l.o.g. that $i=1$ and $B_{p+1}$ is a cross block. Therefore, $B^{-}_{p+1}=b_{p+1}$.
Since $(u,v)$ is a bad arc, $B_{1}$ is a circle block. Consider the path $P_1$ in $\Gamma$ from node $b_{p+1}$ to node $b_1$. 

Let $\Gamma^1$ be the circuit obtained 
by joining the path $P_1$ with the path $\Pi(b_1,u)$ together with arc $(u,v)$, and  $\Pi(v,b_{p+1})$. Note that $\Gamma$ and
$\Gamma^1$ have the same number of row arcs and the same winding number $p$. 
Let us analyze the relationship between $\otimes(\Gamma^1)$ and $\xGam$.

Observe that all nodes that are essential bullets of $\Gamma$, except for $b_{p+1}$, are essential bullets of $\Gamma^1$. Additionally, $v$ is an essential bullet of $\Gamma^1$ which is not an essential bullet of $\Gamma$. 

If we call $B'_i$ with $i\in [s]$ the blocks of $\Gamma^1$, we have that $B_i=B'_i$ for $i\in [s]\setminus \{1,p+1\}$. Moreover, $B'_1\subset B_1$ and $B'_{p+1}$ is a circle block.
Then, $\otimes(\Gamma^1)$ is strictly contained in $\xGam$ 
and the $\Gamma^1$-inequality is stronger than the $\Gamma$-inequality. Hence, the $\Gamma$-inequality does not define a facet of $\Qi{A}$.
\end{myproof}

Now, we can prove:

\begin{lemma}
\label{th:v-not-in-V}
Let $A$ be a circular matrix and $\Gamma$ be a circuit in $F(A)$ with $\xGam \neq \emptyset$, such that the $\Gamma$-inequality is a facet defining inequality for $\Qi{A}$. Let $(u,v)$ be a bad arc with respect to $\Gamma$ with $v \not\in V(\Gamma)$. Then, there is a circuit 
in $F(A)$ defining the same circuit inequality and having less bad arcs than $\Gamma$.
\end{lemma}

\begin{myproof}
We assume w.l.o.g. that $u$ belongs to the circle block $B_1$ and then $v\in (v_{p}, b_{p+1})_n$. From the previous lemma, we know that the block $B_{p+1}$ is not a cross block, as otherwise the $\Gamma$-inequality is not facet defining. Consider the circuit $\Gamma^1$ in $F(A)$ as defined in the previous lemma. 

Since the $\Gamma$-inequality is facet defining, $\otimes(\Gamma^1)=\xGam$ and the $\Gamma^1$-inequality coincides with the $\Gamma$-inequality.  If $\Gamma^1$ has less bad arcs than $\Gamma$, the statement follows. 

Otherwise, since $(u,v)$ is a bad arc with respect to $\Gamma$ but not with respect to $\Gamma_1$, there exists a bad arc $(u^1,v^1)$ with respect to $\Gamma^1$ which is not a bad arc with respect to $\Gamma$.
Since the sets of essential bullets from $\Gamma$ and $\Gamma^1$ differ only in the nodes $b_{p+1}$ and $v$, it follows that $(u^1,v^1)$ must jump over $b_{p+1}$ but not over $v$, i.e., we must have $u^1 \in [v,b_{p+1})_n$ and $v^1\in(b_{2p},b_{2p+1})_n$. 

If $v^1\in V(\Gamma^1)$, by Theorem \ref{th:v-in-V-redundant}, there exists a circuit 
in $F(A)$ with less bad arcs than $\Gamma^1$ defining the same circuit inequality and the statement follows.   

Now, consider the case $v^1 \not\in V(\Gamma^1)$.  
By the previous lemma, $B_{2p+1}$ is not a cross block. Applying iteratively the previous reasoning, either we find a circuit with less bad arcs than $\Gamma$ defining the same circuit inequality, or we obtain that none of the blocks induced by $\Gamma$ is a cross block, contradicting the hypothesis $\xGam \neq \emptyset$.
\end{myproof}

\begin{corollary}
\label{th:v-not-in-V2}
Let $A$ be a circular matrix, $\Gamma$ be a circuit in $F(A)$ with $\xGam \neq \emptyset$ such that the $\Gamma$-inequality is a facet defining inequality for $\Qi{A}$. Let $(u,v)$ be a bad arc with respect to $\Gamma$. Then, there is a circuit $\Gamma'$ without bad arcs such that the $\Gamma'$-inequality coincides with the $\Gamma$-inequality.
\end{corollary}

\begin{myproof}
Due to Theorem \ref{th:v-in-V-redundant} and Lemma~\ref{th:v-not-in-V} it follows that there is a circuit inducing the same inequality as $\Gamma$ and with a less number of bad arcs. Iterating this argument a finite number of times, we prove that there is a circuit $\Gamma'$ without bad arcs which induces the same circuit inequality as $\Gamma$. 
\end{myproof}

Observe that if a circuit $\Gamma$ in $F(A)$ has no crosses then the $\Gamma$-inequality is implied by the rank constraint of $\Qi{A}$. Hence, as a consequence of the previous results, we obtain
for the set covering polyhedron of circular matrices a counterpart of the result obtained by Stauffer \cite{Stauffer05} for the stable set polytope of circular interval graphs. 

\begin{theorem}
\label{th:complete-desc-minors}
Let $A$ be a circular matrix. A complete linear description for the set covering polyhedron $\Qi{A}$ is
given by boolean inequalities, the rank constraint, and $\Gamma$-inequalities with $\Gamma$ a circuit in $F(A)$ without bad arcs.
Moreover, the relevant inequalities for the set covering polyhedron $\Qi{A}$ are minor related row family inequalities induced by circulant minors $\C{s}{p}$ of $A$, with $\gcd(s,p)=1$.
\end{theorem}

\section*{Acknowledgements}
We thank Gianpaolo Oriolo and Gautier Stauffer for earlier discussions that motivated our work on this topic.


\begin{thebibliography}{10}
\bibitem{Aguilera09}  N. Aguilera, \emph{On packing and covering polyhedra of consecutive ones circulant clutters},
Discrete Applied Mathematics (2009), 1343--1356.

\bibitem{ArgiroffoBianchi09} G. Argiroffo, S. Bianchi, \emph{On the set covering polyhedron of circulant matrices}, Discrete Optimization Vol 6-2 (2009), 162--173.

\bibitem{ArgiroffoBianchi10} G. Argiroffo, S. Bianchi, \emph{Row family inequalities for the set covering polyhedron}, Electronic Notes in Discrete Mathematics Vol 36  (2010), 1169--1176.

\bibitem{Bange}Bange,D.W., Barkauskas, A.E., Host, L.H., Slater, P.J., \emph{Generalized domination and efficient domination in graphs}. Discrete Mathematics 159 (1996), 1--11.

 
\bibitem{BartholdiEtAl80} Bartholdi, J. J., Orlin, J. B., and Ratliff, H., \emph{Cyclic scheduling via integer programs with circular ones}. Operations Research Vol 28 (1980), 1074--1085.


\bibitem{Berge62} Berge, C., \emph{Theory of Graphs and its Applications}, Methuen, London (1962).

\bibitem{Bertossi84} Bertossi, A., \emph{Dominating sets for split and bipartite graphs}, Information Processing Letters \textbf{19} (1984), 37--40.

\bibitem{BianchiEtAl14a} S. Bianchi and G. Nasini and P. Tolomei, \emph{Some advances on the set covering polyhedron of circulant matrices}, Discrete Applied Mathematics Vol 166 (2014), 59--60.

\bibitem{BianchiEtAl14b} S. Bianchi and G. Nasini and P. Tolomei, \emph{The minor inequalities in the description of the set covering polyhedron of circulant matrices}, Math.~Meth.~Oper.~Res. Vol 79 (2014), 69--85.

\bibitem{BouchakourEtAl08} M. Bouchakour and T.~M. Contenza and C.~W. Lee and A.~R. Mahjoub, \emph{On the dominating set polytope}, European Journal of Combinatorics Vol 29-3 (2008), 652--661.

\bibitem{Chang98b} Chang, M.-S., \emph{Efficient Algorithms for the Domination Problems on Interval and Circular-Arc Graphs}, SIAM J. Comput., 27(6) (1998), 1671--1694.
 
\bibitem{Chang98} Chang, G., \emph{Algorithmic aspects of domination in graphs}, Handbook of Combinatorial Optimization Vol. 3, D. Du and P Pardalos editors, Kluwer (1998).

\bibitem{ChudnovskySeymour05} M. Chudnovsky and P. Seymour, \emph{The structure of claw-free graphs}, Surveys in combinatorics, London Mathematical Society Lecture Note Series Vol. 327, B.~S. Webb editor,  Cambridge University Press (2005), 153--171.

\bibitem{ChudnovskySeymour08} M. Chudnovsky and P. Seymour, \emph{Claw-free graphs. III. Circular interval graphs}, Journal of Combinatorial Theory, Series B Vol 98-4 (2008), 812--834.

\bibitem{CookEtAl98} W. J. Cook, W.H. Cunningham, W.R. Pulleyblank and A. Schrijver, \emph{Combinatorial Optimization}, John Wiley \& Sons (1998).

\bibitem{CorneilStewart90} Corneil, D. and L. Stewart, \emph{Dominating sets in perfect graphs}, Discrete Mathematics \textbf{86} (1990), 145--164.

\bibitem{Cornuejols94} G. Cornu\'ejols and B. Novick, \emph{Ideal 0 - 1 Matrices}, Journal of Combinatorial Theory
B Vol 60 (1994), 145--157.

\bibitem{EisenbrandEtAl08} F. Eisenbrand and G. Oriolo and G. Stauffer and P. Ventura, \emph{The stable set polytope of quasi-line graphs}, Combinatorica Vol 28-1 (2008), 45--67.

\bibitem{Farber84} Farber, M., \emph{Domination, independent domination and duality in strongly chordal graphs}, Discrete Applied Mathematics \textbf{7} (1984), 115--130.


\bibitem{HaynesEtAl98} Haynes, T., S. Hedetniemi and P. Slater, \emph{Fundamentals of domination in graphs}, Marcel Dekker (1998).

\bibitem{HedetniemiEtAl86} Hedetniemi, S., R. Laskar and J. Pfaff, \emph{A linear algorithm for finding a minimum dominating set in a cactus}, Discrete Applied Mathematics \textbf{13} (1986), 287--292.

\bibitem{KikunoEtAl83} Kikuno, T., N. Yoshida and Y. Kakuda, \emph{A linear algorithm for the domination number of a series-parallel graph}, Discrete Applied Mathematics \textbf{5} (1983), 299--312.

\bibitem{KratschEtAl93} Kratsch, D. and L. Steward, \emph{Domination on comparability graphs}, SIAM. J. on Discrete Mathematics \textbf{6} (3) (1993), 400--417.

\bibitem{Lee} Lee, C.M., Chang, M.S., \emph{Variations of Y-dominating functions on graphs}. Discrete Mathematics 308, (2008), 4185--4204.


\bibitem{Oriolo03} G. Oriolo, \emph{Clique family inequalities for the stable set polytope for quasi-line graphs}, Discrete Applied Mathematics Vol 132-366 (2003), 185--201.

%
\bibitem{PecherWagler06} P{\^e}cher, A. and Wagler ,A.K., \emph{Almost all webs are not rank-perfect}, Mathematical Programming Vol 105, 2 (2006), 311--328.

\bibitem{Sassano} A. Sassano: \emph{On the facial structure of the set covering polytope}, Mathematical Programming Vol 44 (1989), 181--202.

\bibitem{Stauffer05} G. Stauffer, \emph{On the Stable Set Polytope of Claw-free Graphs}, PhD Thesis EPF Lausanne (2005).

\bibitem{TolomeiTorres15} P. Tolomei and L.M. Torres, \emph{Generalized minor inequalities for the set covering polyhedron related to circulant matrices}, Discrete Applied Mathematics (2016) \textbf{210}, 214--222.

\bibitem{Torres15} L.M. Torres, \emph{Minor related row family inequalities for the set covering polyhedron of circulant matrices}, Electronic Notes in Discrete Mathematics Vol 50 (2015), 325--330.

\bibitem{Yanakakis} M. Yannakakis and F. Gavril, \emph{Edge dominating sets in graphs},
SIAM J. Appl. Math. (1980) \textbf{38(3)}, 364--372.


\end{thebibliography}
\end{document}